\documentclass[aos,MSNbibl,seceqn,dvips]{arximspdf}
\usepackage{graphicx}
\doi{10.1214/13-AOS1194} 
\volume{42}
\issue{2}
\pubyear{2014}
\firstpage{592}
\lastpage{624}

\makeatletter

\newcommand{\rrVert}{\Vert}
\newcommand{\rrvert}{\vert}
\newcommand{\llVert}{\Vert}
\newcommand{\llvert}{\vert}

\newtheorem{theorem}{Theorem}
\newproclaim{remark}{Remark}
\newcommand{\argmin}{\mathop{\operatorname{arg\,min}}}
\def\T{\mathrm{T}}
\def\diag{\operatorname{diag}}
\def\qif{\mathrm{QIF}}
\def\pqif{\mathrm{PQIF}}
\def\DZ{\mathrm{DZ}}
\def\DBcalS{{\mathrm{DB}_{\mathrm{S}}}}
\def\DZcalS{{\mathrm{DZ}_{\mathrm{S}}}}
\newcommand{\sgn}{\operatorname{sgn}}
\makeatother

\setattribute{abstract}{width}{278pt}

\begin{document}
\begin{frontmatter}

\title{Estimation and model selection in generalized additive partial
linear models for correlated data with diverging number of covariates\thanksref{TT1}}
\runtitle{Generalized additive partial linear models}

\thankstext{TT1}{The first two authors contribute equally to the paper.}

\begin{aug}
\author[a]{\fnms{Li} \snm{Wang}\thanksref{t1}\ead[label=e1]{lilywang@uga.edu}},
\author[b]{\fnms{Lan} \snm{Xue}\thanksref{t2}\ead[label=e2]{xuel@stat.oregonstate.edu}},
\author[c]{\fnms{Annie} \snm{Qu}\corref{}\thanksref{t3}\ead[label=e3]{anniequ@illinois.edu}}
\and
\author[d]{\fnms{Hua} \snm{Liang}\thanksref{t4}\ead[label=e4]{hliang@gwu.edu}}

\thankstext{t1}{Supported in part by NSF Grants DMS-09-05730,
DMS-11-06816 and DMS-13-09800.}
\thankstext{t2}{Supported in part by NSF Grant DMS-09-06739.}
\thankstext{t3}{Supported in part by NSF Grants DMS-09-06660 and DMS-13-08227.}
\thankstext{t4}{Supported in part by NSF Grant DMS-12-07444 and by
Award Number 11228103, given by National Natural Science Foundation of China.}
\runauthor{Wang, Xue, Qu and Liang}

\affiliation{University of Georgia, Oregon State University,
University of Illinois at Urbana-Champaign and
George Washington University}

\address[a]{L. Wang\\
Department of Statistics\\
University of Georgia\\
Athens, Georgia 30602\\
USA\\
\printead{e1}}

\address[b]{L. Xue\\
Department of Statistics\\
Oregon State University\\
Corvallis, Oregon 97331\\
USA\\
\printead{e2}}

\address[c]{A. Qu\\
Department of Statistics\\
University of Illinois at Urbana-Champaign\\
Champaign, Illinois 61820\\
USA\\
\printead{e3}}

\address[d]{H. Liang\\
Department of Statistics\\
George Washington University\\
Washington, D.C. 20052\\
USA\\
\printead{e4}}
\end{aug}

\received{\smonth{6} \syear{2013}}
\revised{\smonth{11} \syear{2013}}

%
\begin{abstract}
We propose generalized additive partial linear models for complex data
which allow one to capture nonlinear patterns of some covariates, in
the presence of linear components. The proposed method improves
estimation efficiency and increases statistical power for correlated
data through incorporating the correlation information. A~unique
feature of the proposed method is its capability of handling model
selection in cases where it is difficult to specify the likelihood
function. We derive the quadratic inference function-based estimators
for the linear coefficients and the nonparametric functions when the
dimension of covariates diverges, and establish asymptotic normality
for the linear coefficient estimators and the rates of convergence for
the nonparametric functions estimators for both finite and
high-dimensional cases. The proposed method and theoretical development
are quite challenging since the numbers of linear covariates and
nonlinear components both increase as the sample size increases. We
also propose a doubly penalized procedure for variable selection which
can simultaneously identify nonzero linear and nonparametric
components, and which has an asymptotic oracle property. Extensive
Monte Carlo studies have been conducted and show that the proposed
procedure works effectively even with moderate sample sizes. A
pharmacokinetics study on renal cancer data is illustrated using the
proposed method.
\end{abstract}

%
\end{frontmatter}

\section{Introduction}\label{sec:intr}

We encounter longitudinal data in many social and health studies where
observations from clustered data are measured over time, and can often
be discrete, such as binary or count data. Generalized additive partial
linear models (GAPLM) are developed to model partial linear additive
components while the remaining components are modeled nonparametrically
\cite{har:mul:spe:04} to combine the strengths of both the GPLM and
the GAM for interpretability and flexibility.

Efficient estimation of linear and nonparametric function components is
quite challenging even for cross-sectional data. To solve the ``curse
of dimensionality'' problem in computing, \cite{Wood:JASA:04} suggested
a penalized regression splines approach to utilize the practical
benefits of smoothing spline methods and the computational advantages
of local scoring backfitting \cite{buj:has:tib:89}. In addition,
\cite
{wang:liu:liang:carroll:11} applied polynomial splines to approximate
the nonparametric components, and estimated coefficients through an
efficient one-step procedure of maximizing the quasi-likelihood
function. This can reduce computational costs significantly compared to
the local scoring backfitting and marginal integration approaches.
Another advantage of the polynomial spline approach is that it can
formulate a penalized function for variable selection purposes, which
cannot be easily implemented through other iterative methods.

However, \cite{wang:liu:liang:carroll:11}'s approach is valid only for
independent data and the case with a fixed number of covariates for
linear component model selection. In this paper, we develop a general
framework for estimation and variable selection using the GAPLM. The
proposed method can handle correlated categorical responses in addition
to continuous ones, and allows both the number of covariates for linear
and nonlinear terms to diverge as the sample size increases. Note that
the theoretical development for model selection and estimation for
diverging number of covariates in nonlinear components are completely
different from the setting with finite dimension of covariates \cite
{xue:qu:zhou:10}.

The GAPLM can be highly computationally intensive as it introduces
high-dimensional nuisance parameters associated with nonparametric
forms. Incorporating correlation structure brings additional challenges
to modeling and estimation due to the additional correlation parameters
involved. The extension of the GAPLM for correlated data imposes more
challenges computationally and theoretically. However, it is well known
that ignoring correlation could lead to inefficient estimation and
diminish statistical power in hypothesis testing and the selection of
correct models. Moreover, \cite{wang:naisyin:03} and \cite
{Zhu:Fung:He:08} indicate that in nonparametric settings ignoring the
correlation could also result in biased estimation since the selection
process is rather sensitive to small departures from the true
correlation structure, and likely to cause overfitting of the
nonparametric estimator to compensate for the overall bias. These
problems could be more critical for the GAPLM since in contrast to the
parametric setting, the true model here might be more difficult to
verify. The proposed polynomial spline approach can efficiently take
the within-cluster correlation into account because of its nonlocal
behavior in longitudinal data \cite{welsh:lin:carroll:02}. This is
substantially different from the kernel smoothing method, where only
local data points are used in the estimation and, therefore, it cannot
incorporate correlation structure efficiently.

We propose variable selection and estimation simultaneously based on
the penalized quadratic inference function for correlated data when the
dimension of covariates in GAPLM increases as the sample size. The
quadratic inference function~(QIF) \cite{qu:lindsay:li:00} utilizes
within-cluster correlation into account without specifying the
likelihood function, and is less sensitive to the misspecification of
working correlation matrices compared to the generalized estimating
equation (GEE) method \cite{liang:zeger:86}, in general. In addition,
we perform variable selection for the marginal GAPLM to identify
important variables, which is crucial to obtain efficient estimators
for the nonzero components. We show that the proposed model selection
for both parametric and nonparametric terms is consistent, the
estimators of the nonzero linear coefficients are asymptotically
normal, and the estimators of the nonzero nonparametric functions are
$L_{2}$-norm consistent with the optimal rate of convergence if the
dimension of nonparametric components is finite. However, the
asymptotic properties on the rate of convergence are no longer the same
as in \cite{wang:liu:liang:carroll:11} when the dimensions of
covariates for parametric and nonparametric components both diverge as
the sample size increases.

The semiparametric model containing both linear and nonparametric
functions makes the estimation and model selection very different from
the generalized additive model \cite{xue:qu:zhou:10}, which involves
only nonparametric components. The establishment of the asymptotic
normal distribution of the estimators for the parametric terms is quite
challenging given that the number of covariates for both parametric and
nonparametric terms diverge, and the convergence rate for the
nonparametric component estimators is slower than $\sqrt{n}$. Another
difficulty here is that the covariates in the parametric components and
those in the nonparametric components could be dependent, in addition
to dependent errors for repeated measurements, so traditional
nonparametric tools such as the backfitting algorithm \mbox{\cite
{buj:has:tib:89}} cannot be applied here. In contrast, the proposed
spline-based approach allows one to incorporate correlation effectively
even when the number of covariates diverges.

In addition, the required techniques using the penalized quadratic
distance function for the diverging numbers of linear and nonlinear
covariates setting are very different from existing approaches such as
the penalized least-squares approach for a finite dimension setting
\cite{xue:09,ma:song:wang:11,wang:liu:liang:carroll:11}; the
generalized linear model selection approach for the parametric term
only with diverging number of covariates \mbox{\cite{cho:qu:13}}; or
the GAPLM for a finite number of nonparametric functions \mbox{\cite
{lian:liang:wang:13}}, which does not perform model selection for the
nonparametric term. This motivates us to develop new theoretical tools
to derive large sample properties for linear and nonparametric
components estimation and model selection to incorporate the dependent
nature of the data for handling diverging numbers of covariates.

We organize the paper as follows. Section~\ref{sec:modelsec} presents
the model framework, describes estimation procedures, and establishes
asymptotic properties of the GAPLM for correlated data. Section~\ref{sec:pQIF} proposes a penalized QIF method for simultaneous estimation
and variable selection when the dimension of covariates increases as
the sample size. The theoretical properties on model selection
consistency and rate of convergence for the nonparametric estimators
are developed, in addition to algorithm implementation and tuning
parameter selection. Sections~\ref{sec:simu} and \ref{sec:real}
illustrate the performance of the proposed method through simulation
studies and a pharmacokinetics study on renal cancer patients,
respectively. We provide concluding remarks and discussion in Section~\ref{sec:dics}. The proofs of the theorems along with technical lemmas
are provided in the \hyperref[app]{Appendix} and supplementary material~\cite
{wang:xue:qu:liang:13}.

\section{Estimation procedures and theoretical results}\label{sec:modelsec}

\subsection{The GAPLM for correlated data}

For the clustered data, let $Y_{it}$ be a response variable, $%
\mathbf{X}_{it}=( X_{it}^{(1)},\ldots,X_{it}^{(d_{x})})^{\T}$ and
$\mathbf{Z%
}_{it}=(1,Z_{it}^{(1)},\ldots,Z_{it}^{(d_{z}-1)})^{\T}$ be the $d_{x}$-vector
and $d_{z}$-vector of covariates corresponding to the nonparametric and
parametric components, respectively, where $t$ is the $t$th ($t=1,\ldots
,T_{i}$) observation for the $i$th ($i = 1, \ldots, n$) cluster. Further
denote $\mathbf{Y}_{i}= ( Y_{i1},\ldots,Y_{iT_{i}} )
^{\T
}$, $%
\mathbf{X}_{i}= ( \mathbf{X}_{i1},\ldots,\mathbf
{X}_{iT_{i}} ) ^{\T%
}$, and $\mathbf{Z}_{i}= ( \mathbf{Z}_{i1},\ldots
,\mathbf{Z}%
_{iT_{i}} ) ^{\T}$. For presentation simplicity, we
assume each cluster has the same size with $T_{i}=T<\infty$. The
procedure for
data with unequal cluster sizes can be adjusted following the same
method of
\cite{xue:qu:zhou:10}.

One of the advantages of marginal approaches is that we only need to specify
the first two moments by $E ( Y_{it}|\mathbf
{X}_{it},\mathbf
{Z}%
_{it} ) =\mu_{it}$, and $\mathrm{Var}(Y_{it}|\mathbf
{X}_{it},\mathbf{Z}%
_{it})=\phi V ( \mu_{it} ) $, where $\phi$ is a scale parameter
and $V ( \cdot ) $ is a known variance function. Here, the marginal
mean $\mu_{it}$ associates with the covariates through the known link
function $g ( \cdot )$ such that
%
\begin{equation}
\eta_{it}=g ( \mu_{it} ) =\sum
_{l=1}^{d_{x}}\alpha _{l} \bigl(
X_{it}^{(l)} \bigr) +\mathbf{Z}_{it}^{\T}
\bolds{\beta}, \label{mgaplm}
\end{equation}
where $\bolds{\beta}$ is $d_{z}$-vector of unknown parameters,
and $%
\{ \alpha_{l} ( \cdot )  \} _{l=1}^{d_{x}}$
are unknown
smooth functions. Model (\ref{mgaplm}) is called the generalized
additive partial linear model (GAPLM), where $\mathbf{Z}_{it}^{\T
}\bolds{%
\beta}$ are the parametric components, and $\sum_{l=1}^{d_{x}}\alpha
_{l}(X_{it}^{(l)})$ are the nonparametric components. Here, the
mean of $Y_{it}$ depends only on the covariate vector for the %
$t$th observation \cite{pepe:anderson}, that is, $E(Y_{it}|%
\mathbf{X}_{i}, \mathbf{Z}_{i})= E(Y_{it}|\mathbf{X}_{it},
\mathbf{Z}_{it})$%
. In addition, without loss of generality, we assume that each
covariate $ \{X^{(l)} \}_{l=1}^{d}$ can be rescaled into
$[0,1]$; and each $\alpha
_{l} ( \cdot )$ is centered with $\int_{0}^{1}\alpha
_{l} (
x ) \,dx=0$ to make model~(\ref{mgaplm}) identifiable.

\subsection{Spline approximation}\label{subsec:splineQIF}

We approximate smooth functions $ \{ \alpha_{l} ( \cdot
)
\} _{l=1}^{d_{x}}$ in~(\ref{mgaplm}) by polynomial splines
for their simplicity in computation, and they often provide a good
approximation of smooth functions with a limited number of knots. For
example, for each $1\leq l\leq d_{x}$, let $\upsilon_{l}$ be a
partition of $[0,1]$, with $N_{n}$ interior knots
$\upsilon_{l}= \{ 0=\upsilon_{l,0}<\upsilon_{l,1}<\cdots
<\upsilon
_{l,N_{n}}<\upsilon_{l,N_{n}+1}=1 \}$.

The polynomial splines of order $p+1$ are functions with $p$-degree (or
less) of polynomials on intervals $[\upsilon_{l,i},\upsilon
_{l,i+1}),i=0,\ldots,N_{n}-1$, and $ [ \upsilon
_{l,N_{n}},\upsilon
_{l,N_{n}+1} ] $, and have $p-1$ continuous derivatives globally. Let
$\varphi_{l}=\varphi^{p}( [ 0,1 ],\upsilon_{l})$ be the
space of such polynomial splines, and $\varphi_{l}^{0}= \{ s\in
\varphi_{l}\dvtx \int_{0}^{1}s ( x ) \,dx=0 \}$. This ensures
that the
spline functions are centered.

Let $ \{ B_{lj} ( \cdot )  \} _{j=1}^{J_{n}}$ be a
set of
spline bases of $\varphi_{l}^{0}$ with the dimension of
$J_{n}=N_{n}+p$. We
approximate the nonparametric component ${\alpha}${$_{l} ( \cdot
) $} by a polynomial spline, that is $\alpha_l(\cdot) \approx
s_{l} ( \cdot ) =\sum_{j=1}^{J_{n}}\gamma_{lj}B_{lj} (
\cdot
) $, with a set of coefficients $\bolds{\gamma
}_{l}= (
\gamma
_{l1},\ldots,\gamma_{lJ_{n}} ) ^{\T}$. Accordingly, $\eta
_{it}$ is
approximated by
\[
\eta_{it}(\bolds{\beta},\bolds{\gamma})=\sum
_{l=1}^{d_{x}}%
\sum
_{j=1}^{J_{n}}\gamma_{lj}B_{lj}
\bigl( X_{it}^{(l)} \bigr) +\mathbf{Z}%
_{it}^{\T}\bolds{\beta},
\]
where $\bolds{\gamma}= ( \bolds{\gamma}_{1}^{\T
},\ldots
,\bolds{%
\gamma}_{d_{x}}^{\T} ) ^{\T}$. Therefore, the mean function $%
\mu_{it}$ in (\ref{mgaplm}) can be approximated by
\[
\mu_{it}\approx\mu_{it}(\bolds{\beta},\bolds{
\gamma })=g^{-1} \Biggl\{ \sum_{l=1}^{d_{x}}
\sum_{j=1}^{J_{n}}\gamma_{lj}B_{lj}
\bigl( X_{it}^{(l)} \bigr) +\mathbf{Z}_{it}^{\T}
\bolds{\beta } \Biggr\}.
\]
We denote $\bolds{\mu}_{i} (\bolds{\beta
},\bolds{\gamma} )
= \{ \mu_{i1}(\bolds{\beta},\bolds{\gamma}),\ldots
,\mu_{iT}(%
\bolds{\beta},\bolds{\gamma}) \}^{\T}$ in matrix
notation. To
incorporate the within-cluster correlation, we apply the QIF to estimate
$\bolds{\beta}$ and $\bolds{\gamma}$ for the
parametric and nonparametric parts, respectively.

\subsection{Quadratic inference functions}\label{subsec:QIFp}

To estimate $\bolds{\beta}$ and $\bolds{\gamma}$, one may use
the GEE method \cite{lian:liang:wang:13}, that is, using a working
correlation matrix $R$ which depends on fewer nuisance parameters. The
estimates of regression parameters are consistent even when $R$ is
misspecified. However, one has to find a consistent estimator of $R$ to
obtain an efficient estimator of $\bolds{\beta}$. The QIF approach
\cite{qu:lindsay:li:00} considers the approximation of $\mathbf
{R}^{-1}$ with a linear combination of basis matrices of form
$\mathbf{R}^{-1}\approx a_{1}\mathbf{M}_{1}+\cdots
+a_{K}\mathbf{M}_{K}$. For example, if $\mathbf{R}$ has an
exchangeable structure with correlation $\rho$, then $\mathbf
{R}^{-1}$ can be represented as $a_{1} \mathbf{I}+a_{2}\mathbf
{M}_{2}$ with $\mathbf{I}$ being the identity matrix and
$\mathbf{M}_{2}$ being a matrix with $0$ on the diagonal and $1$
off the diagonal. The corresponding coefficients are $a_{1}=- \{
( T-2 ) \rho+1 \} /k_{1}$, and $a_{2}=\rho/k_{1}$,
where $k_{1}=(T-1)\rho^{2}-(n-2)\rho-1$ and $T$ is the dimension of
$\mathbf{R}$. The basis matrices are also available to approximate
$\mathbf{R}^{-1}$ of other structure such as, AR-1 and the block
diagonal correlation structures. If the candidate basis matrices
represent a sufficiently rich class for the true structure, \cite
{zhou:qu:12} show that the correlation structure can be selected
consistently by minimizing the penalized difference between two
estimating functions generated from the empirical correlation
information and the model-based approximation, respectively. The
penalization on the basis matrices ensures that an optimal number of
basis matrices $K$ will be selected to capture correlation information,
yet not be burdened by too many moment conditions.

The quadratic inference function is established under the same
principle as the generalized method of moments \mbox{\cite
{hansen:82}}, and is shown to be the most efficient among estimators
given the same class of estimating functions as the asymptotic variance
reaches the minimum in the sense of Loewner ordering. This is
especially useful under misspecified working correlation structures,
since the true correlation structure is seldom known. For example, the
QIF estimator is shown to be more efficient than the GEE estimator for
diverging number of covariates under the generalized linear model
framework \mbox{\cite{cho:qu:13}}. Another advantage of the QIF is
that the estimation of the linear
coefficients $a_{i}$'s is not required. In nonparametric modeling with
diverging number of covariates, it is even more beneficial if we can
avoid estimating the nuisance parameters associated with the
correlations, since we are dealing with high-dimensional parameters
involved in nonparametric components.

\subsection{Estimation procedure}

$\!\!$For any $\mathbf{x}\in R^{d_x}$, $\mathbf{z}\in R^{d_z}$, let $
\mathbf{B}^{\T}(\mathbf{x})=(B_{11}(x_{1}),\break   \ldots
,B_{1J_{n}}(x_{1}), \ldots, {B_{d_{x}1}}(x_{d_{x}}),\ldots
, B_{d_{x}J_{n}}(x_{d_{x}}))$, $\mathbf{D}^{\T}(\mathbf
{x},\mathbf{z})=(\mathbf{z}^{\T},\mathbf{B}^{\T
}(\mathbf{x}))$
be vectors of dimensions $d_{x}J_{n}$ and $d_{x}J_{n}+d_{z}$,
respectively. In addition, we denote matrices $\mathbf
{B}_{i}=
\{(\mathbf{B}(\mathbf{X}_{i1}),\ldots,\mathbf
{B}(\mathbf
{X}_{iT}))^{\T} \} _{T\times d_{x}J_{n}}$, $\mathbf
{D}_{i}= \{  (\mathbf{D}(\mathbf
{X}_{i1},\mathbf
{Z}_{i1}), \ldots,\break \mathbf
{D}(\mathbf{X}_{iT}, \mathbf{Z}_{iT}) )^{\T} \}
_{T\times (
d_{x}J_{n}+d_{z} )}$.

For $\bolds{\theta}= ( \bolds{\beta}^{\T
},\bolds{\gamma}^{\T} ) ^{\T}$,
we define $K(d_xJ_n+d_z)$-dim extended scores to incorporate correlation
for correlated data as follows:
%
\begin{equation}
\mathbf{g}_{i}(\bolds{\theta})= \pmatrix{ \mathbf{D}_{i}^{\T}\bolds{
\Delta}_{i}\mathbf {A}_{i}^{-1/2}
\mathbf{M}_{1}%
\mathbf{A}_{i}^{-1/2}
\bigl\{ \mathbf{Y}_{i}-\bolds{\mu }_{i}(
\bolds{\theta }) \bigr\}
\vspace*{2pt}\cr
\vdots
\vspace*{2pt}\cr
\mathbf{D}_{i}^{\T}\bolds{\Delta}_{i}
\mathbf {A}_{i}^{-1/2}\mathbf{M}_{K}%
\mathbf{A}_{i}^{-1/2} \bigl\{ \mathbf{Y}_{i}-
\bolds{\mu }_{i}(\bolds{\theta }) \bigr\}}, \label{DEF:githeta}
\end{equation}
where $\bolds{\Delta}_{i}=\diag \{ \dot{\mu}_{i1},\ldots
,\dot
{\mu}%
_{i} \} $ {and} $\dot{\mu}_{it}$ is the first order derivative
of $
g^{-1}$ evaluated at $\mathbf{B}^{\T} ( \mathbf
{X}_{it} ) \bolds{%
\gamma}+\mathbf{Z}_{it}^{\T}\bolds{\beta}$; and
$\mathbf
{A}_{i}=\diag%
\{ V ( \mu_{i1} ),\ldots,V ( \mu_{i} )
\} $. We define the sample mean and sample variance of the moment
conditions as
%
\begin{equation}
\mathbf{G}_{n}(\bolds{\theta})=\frac{1}{n}\sum
_{i=1}^{n}\mathbf{g}_{i}(%
\bolds{\theta}),\qquad \mathbf{C}_{n}(\bolds{\theta })=
\frac
{1}{n}%
\sum_{i=1}^{n}
\mathbf{g}_{i}(\bolds{\theta})\mathbf
{g}_{i}^{\T}(\bolds{%
\theta}). \label{DEF:GnCn}
\end{equation}
If we set $\mathbf{G}_{n}(\bolds{\theta})=0$ as our
estimating equations,
there are more equations than the number of unknown parameters, and the
parameters are over-identified. The QIF approach estimates
$\alpha_{l} ( \cdot ) $ and $\bolds{\beta}$ by making
$\mathbf{G}_{n}$
as close to zero as possible, in the sense of minimizing the QIF $Q_{n}(
\bolds{\theta})$, that is,
%
\begin{eqnarray} \label{EQ:thetahat_qif}
\widehat{\bolds{\theta}}^{\qif} &=& \bigl(\bigl(\widehat {\bolds{\beta}}^{\qif}\bigr)^{\T},\bigl(%
\widehat{
\bolds{\gamma}}^{\qif}\bigr)^{\T}\bigr)^{\T}
\nonumber
\\[-8pt]
\\[-8pt]
\nonumber
& =&
\argmin _{\bolds{%
\theta}}Q_{n} ( \bolds{\theta} ) = \argmin
_{\bolds{%
\theta}} \bigl\{ n\mathbf{G}_{n}^{\T}(
\bolds{\theta })\mathbf{C}%
_{n}^{-1} (
\bolds{\theta} ) \mathbf {G}_{n}(\bolds{\theta}%
) \bigr\}.
\end{eqnarray}

Consequently, for any $\mathbf{x}\in{}[0,1]^{d_{x}}$ and
$l=1,\ldots,d_{x}$, {the estimators
of the nonparametric components in (\ref{mgaplm}) are provided as}
%
\begin{equation}
\widehat{\alpha}_{l}^{\qif} \bigl( x^{(l)} \bigr) =
\sum_{j=1}^{J_{n}}%
\widehat{
\gamma}_{lj}^{\qif}B_{lj} \bigl( x^{ ( l )
}
\bigr)\quad \mbox{and}\quad \widehat{\alpha}^{\qif} ( \mathbf{x}%
) =\sum_{l=1}^{d_{x}}\widehat{
\alpha}_{l}^{\qif} \bigl( x^{ (
l ) } \bigr).
\label{EQ:alphahat}
\end{equation}

The advantages of the spline basis approach lie not only in its
computation efficiency, but also in the ease of implementation. Using
the spline basis approximation, we can easily convert a problem with
infinite-dimensional parameters to one with a finite number of
parameters \cite{huang:zhang:zhou:07}. In the following Theorem~\ref{THM:nonparametric}, we
also show that the proposed estimators of the nonparametric components
using polynomial spline achieve the optimal rate of convergence. This
result is useful for providing an initial consistent estimator for
later development in simultaneous variable selection and estimation for
both parametric and nonparametric functions.

\subsection{Asymptotic properties}\label{subsec:asym}

We establish the asymptotic properties of the QIF estimators, summarize the
main results in the following theorems and provide detailed proofs in
the \hyperref[app]{Appendices}. Note that the asymptotic results still hold for unequal
cluster size data.

In the following, denote the true nonparametric components by $\alpha
_{0,l}$%
, $1\leq l\leq d_{x}$ and the true parameters for the parametric components
by $\bolds{\beta}_{0}$. Let $\mu_{0,it}$ be the true marginal
means. In
addition, let $\bolds{\mu}_{0,i}= ( \mu_{0,i1},\ldots,\mu
_{0,iT} ) ^{\T}$ and $\mathbf{e}_{i}=\mathbf
{Y}_{i}-\bolds{\mu}_{0,i}$. Let $\bolds{\Gamma}_{0,i}^{(k)}=\bolds{\Delta
}_{0,i}\mathbf{V}%
_{0,i}^{(k)}\bolds{\Delta}_{0,i}$, where $\mathbf
{V}_{0,i}^{(k)}=\mathbf{A}_{0,i}^{-1/2}\mathbf
{M}_{k}\mathbf
{%
A}_{0,i}^{-1/2}$ and $\bolds{\Delta}_{0,i}$, $\mathbf{A}_{0,i}$ are
evaluated at $\bolds{\mu}_{0,i}$. Similarly, define $\bolds{\mu}_{0}$, $%
\mathbf{e}$, $\bolds{\Gamma}_{0}^{(k)}$ as the generic
versions of $%
\bolds{\mu}_{0,i}$, $\mathbf{e}_{i}$ and $\bolds{\Gamma
}_{0,i}^{(k)}$,
respectively, for $ (\mathbf{Y},\mathbf{X}, \mathbf
{Z} )$.
Let $d_{n}=d_xJ_n+d_z$, and $\rho_n=((1-\delta)/2)^{(d_x-1)/2}$, for
some constant $\delta\in(0,1)$. Further, we denote $a\asymp b$, if
there exist constants $c\geq c^{\ast}>0$ such that $c^{\ast}b\leq
a\leq cb$.

\begin{theorem}
\label{THM:nonparametric} Under conditions \textup{(C1)--(C3)}, \textup{(C5)--(C8)} in
Appendix \ref{subsec:assu}, if $d_x/ \log(n)$ $\rightarrow
0$, $n^{-1/4}d_z\rightarrow0$, $J_{n}\asymp n^{b}$, for some
$1/(4r)\leq b
< 1/4$ with the smoothing parameter $r>1$ defined in condition \textup{(C1)},
the estimators $\widehat{\alpha}_{l}^{\qif}( x^{(l)})$, $1\leq l\leq d_x$,
defined in (\ref{EQ:alphahat}) satisfy
\[
{\frac{1}{n}\sum_{l=1}^{d_x}\sum
_{i=1}^{n}\sum
_{t=1}^{T}%
\bigl\{ \widehat{
\alpha}_{l}^{\qif}\bigl( x_{it}^{(l)}\bigr)
-\alpha_{0,l}\bigl( x_{it}^{(l)}\bigr) \bigr\}
^{2}=O_{P}\bigl(n^{-1}d_n+J_{n}^{-2r}d_x
\bigr)},
\]
where $r$ determines the smoothness of the nonparametric functions.
In particular, if $J_{n}\asymp n^{1/(2r+1)}$ and $d_{z}=O(J_{n}d_{x})$, then
\[
{ \frac{1}{n}\sum_{l=1}^{d_x}\sum
_{i=1}^{n}\sum
_{t=1}^{T}%
\bigl\{ \widehat{
\alpha}_{l}^{\qif}\bigl( x_{it}^{(l)}\bigr)
-\alpha_{0,l}\bigl( x_{it}^{(l)}\bigr) \bigr\}
^{2}=O_{P} \bigl(n^{-2r/ ( 2r+1
)}d_x
\bigr)}.
\]
\end{theorem}

\begin{remark}
{Note that $d_{n}=d_xJ_n+d_z$, so if the number of nonparametric
functions, $d_x$, is finite,
and $J_n\asymp n^{1/(2r+1)}$, then we obtain an optimal convergence
rate $n^{-2r/ ( 2r+1 )}$.
In addition, for a cluster size of 1, this reduces to a special case
where the responses
are independent, and is the same as in \cite{huang:98} and \cite
{xue:09} for independent data.}
\end{remark}

Next, we establish the asymptotic normal distribution for
the parametric estimator. We denote $\mathbf{g}_{0,i}=(\mathbf
{g}_{0,i1}^{\T%
},\ldots,\mathbf{g}_{0,iK}^{\T})^{\T}$ with\vspace*{1pt} $\mathbf
{g}_{0,ik}=\mathbf{D}_{i}^{%
\T}\bolds{\Delta}_{0,i}\mathbf{V}_{0,i}^{(k)}\mathbf
{e}_{i}$, the value of
$\mathbf{g}_{i}$ in (\ref{DEF:githeta}) at $\bolds{\mu
}_{i}=\bolds{\mu}%
_{0,i}$. Similarly, let
%
\begin{equation}\qquad
\mathbf{G}_{n}^{0}=\frac{1}{n}\sum
_{i=1}^{n}\mathbf {g}_{0,i},\qquad
\mathbf{C}%
_{n}^{0}=\frac{1}{n}\sum
_{i=1}^{n}\mathbf{g}_{0,i}
\mathbf {g}_{0,i}^{\T},\qquad \mathbf{Q}_{n}^{0}=n
\bigl( \mathbf{G}_{n}^{0} \bigr) ^{\T
} \bigl(
\mathbf{C}_{n}^{0} \bigr) ^{-1}
\mathbf{G}_{n}^{0} \label
{DEF:Gn0_Cn0_Qn0}
\end{equation}
be the corresponding values of $\mathbf{G}_{n}$, $\mathbf
{C}_{n}$ and $%
\mathbf{Q}_{n}$ defined in (\ref{DEF:GnCn}) and (\ref
{EQ:thetahat_qif}) at \mbox{$%
\bolds{\mu}_{i}=\bolds{\mu}_{0,i}$}. Next, denote
$\widehat
{\mathbf{Z}}_{i}=\mathbf{Z}_{i}-\operatorname{Proj}_{\Gamma
_{n}}\mathbf{Z}%
_{i}$,
where $\operatorname{Proj}_{\Gamma_{n}}$ is the projection onto the empirically centered
additive spline space.
See (S.17) for the exact formula of~$\widehat{\mathbf{Z}}_{i}$.
Further denote
%
\begin{eqnarray}
\widehat{\mathbf{J}}_{\DZ}^{(k)} &=&\frac{1}{n}\sum
_{i=1}^{n}\mathbf{D}%
_{i}^{\T}\bolds{\Gamma}_{0,i}^{ ( k ) }
\widehat {\mathbf{Z}}%
_{i},\qquad \widehat{
\mathbf{J}}_{\DZ}= \bigl\{ \bigl(\widehat {\mathbf
{J}}_{\DZ%
}^{(1)}\bigr)^{\T},\ldots,\bigl(\widehat{
\mathbf{J}}_{\DZ}^{(K)}\bigr)^{\T
} \bigr\}
^{\T}, \label{DEF:JhatDZ}
\\
\mathbf{W}_{i}^{(k)} &=&\mathbf{D}_{i}^{\T}
\bolds{\Delta }_{0,i}\mathbf{V}%
_{0,i}^{ ( k ) }\bolds{\Sigma}_{i}^{1/2},\qquad
\mathbf{W}%
_{i}= \bigl\{ \bigl(\mathbf{W}_{i}^{(1)}
\bigr)^{\T},\ldots,\bigl(\mathbf {W}_{i}^{(K)}
\bigr)^{\T%
} \bigr\} ^{\T}. \label{DEF:Wi}
\end{eqnarray}
In what follows, $\mathbf{A}^{\otimes2}$ and $\mathbf
{A}_{\mathbf{B}}^{\otimes2}$ stand for $\mathbf{AA}^{\T}$ and
$\mathbf{ABA}^{\T}$ for any matrix/vector $\mathbf{A}$ and
square matrix $\mathbf{B}$, respectively.

\begin{theorem}
\label{THM:parametric} Assume that conditions \textup{(C1)--(C3)}, \textup{(C5)--(C9)} in
Appendix \ref{subsec:assu} are satisfied,
if $d_x/ \log(n)\rightarrow0$, $n^{-1/5}d_z\rightarrow0$,
and $J_{n}\asymp n^{b}$, for some $1/(2r+1)\leq b < 1/5$, where the
smoothing parameter $r>2$,
then the estimator $\widehat{\bolds{\beta}}^{\qif}$ of
$\bolds{\beta}%
_{0}$ is consistent and
$ 
\sqrt{n}\mathbf{A}_n\bolds{\Sigma}_{n}^{-1/2}(\widehat
{\bolds{\beta}}^{%
\qif}-\bolds{\beta}_{0})\mathop{\rightarrow}^{D} N(0,\bolds{\Sigma}_{A})$, 
where $\mathbf{A}_n$ is any $q \times d_z$ matrix with a finite
$q$ such
that $\mathbf{A}_{n}^{\otimes2}$ converges to a $q\times q$ nonnegative
symmetric $\bolds{\Sigma}_{A}$, and $\bolds{\Sigma
}_{n}=\widehat{\bolds{%
\Psi}}_{n}^{-1}\widehat{\bolds{\Omega}}_{n}\widehat
{\bolds{\Psi}}%
_{n}^{-1}$ with
%
\begin{eqnarray}
\widehat{\bolds{\Psi}}_{n} =\widehat{\mathbf{J}}_{\DZ
}^{\T
}
\bigl( \mathbf{%
C}_{n}^{0} \bigr)
^{-1}\widehat{\mathbf{J}}_{\DZ} \quad\mbox{and}\quad
\widehat{%
\bolds{\Omega}}_{n} =\frac{1}{n}\sum
_{i=1}^{n} \bigl\{ \widehat {
\mathbf{J}}%
_{\DZ}^{\T} \bigl(
\mathbf{C}_{n}^{0} \bigr) ^{-1}\mathbf
{W}_{i} \bigr\} ^{\otimes2}. \label{DEF:OMIGAn}
\end{eqnarray}
\end{theorem}

To establish the asymptotic properties of the QIF estimators for
diverging number of covariates, a crucial step is to obtain the upper
and lower bounds of the eigenvalues of the matrix $\mathbf
{C}_{n}^{-1} ( \bolds{\theta} )$ in (\ref{DEF:GnCn}) and
(\ref{EQ:thetahat_qif}). Note that $\mathbf{C}_{n} (
\bolds{\theta} )$ is a random matrix with increasing dimension of
linear and nonlinear components as $n$ increases. The derivation of its
bounds relies heavily on Lemma~1 of \cite{stone:85}; see \cite
{xue:qu:zhou:10,wang:liu:liang:carroll:11}. When $d_x$ is finite, the
term $\rho_n=((1-\delta)/2)^{(d_x-1)/2}$ in Lemma~1 of~\cite{stone:85}
is a constant, which makes the derivation of the bounds relatively
easy. However, this is no longer true in the diverging case since $\rho
_n$ goes to zero as $d_x$ goes to infinity, and it requires special
techniques for asymptotic derivations. Another major difficulty in the
derivation of Theorem~\ref{THM:parametric} is to resolve the dependence between
$\mathbf{X}$ and $\mathbf{Z}$ in addition to establishing the
convergence results for the first- and second-order partial derivatives
of the quadratic inference
function, which could be infinite-dimensional.

\section{Penalized QIF for marginal GAPLM}\label{sec:pQIF}

In this section, we define predictor variables $X_{l}$ and $Z_k$ as
redundant in model (\ref{mgaplm}), if and only if $\alpha_{l}(X_{l})=0$
and \mbox{$\beta_k=0$}. Suppose there is only an unknown subset of predictor
variables which is relevant in model (\ref{mgaplm}) with nonzero
components, we are interested in identifying such subsets of relevant
predictors consistently while estimating the nonzero parameters and
functions in ({\ref{mgaplm}}) simultaneously.

\subsection{Model selection}\label{subsec:pQIF}

To perform model selection for the GAPLM, we propose the penalized quadratic
inference function in (\ref{EQ:thetahat_qif}) which shrinks small components
of estimated functions to zero. Through consistent model selection, we
are able to improve
the efficiency of estimators for the nonzero components since the
correlation within clusters is taken into account. We define the penalized
QIF (PQIF) estimator as
\begin{eqnarray*}
&& \bigl( \bigl(\widehat{\bolds{\beta}}^{\pqif}\bigr)^{\T},
\bigl(\widehat {\bolds{\gamma}}^{%
\pqif}\bigr)^{\T} \bigr)
^{\T}
\\
&&\qquad =\argmin_{\bolds{\beta},\bolds{\gamma}%
} \Biggl\{ Q_{n} ( \bolds{\beta},
\bolds{\gamma} ) +n\sum_{l=1}^{d_{x}}p_{\lambda_{1,n}}\bigl(
\Vert\bolds{\gamma }_{l}\Vert_{%
\mathbf{K}_{l}}\bigr)+n\sum
_{l=1}^{d_{z}}p_{\lambda_{2,n}}\bigl(|\beta _{l}|\bigr)
\Biggr\}, \label{DEF:pqif}
\end{eqnarray*}
where $p_{\lambda_{\bullet, n}} ( \cdot ) $ are given penalty
functions of tuning parameters $\lambda_{\bullet, n}$, and $\Vert
\bolds{%
\gamma}_{l}\Vert_{\mathbf{K}_{l}}^{2}=\bolds{\gamma
}_{l}^{\T
}\mathbf{K}%
_{l}\bolds{\gamma}_{l}$, in which $\mathbf{K}_{l}=\frac
{1}{n}\sum_{i=1}^{n}%
\frac{1}{T}\sum_{t=1}^{T}\mathbf{B}_{l} ( \mathbf
{X}_{it}^{ (
l ) } ) \mathbf{B}_{l}^{\T} ( \mathbf
{X}_{it}^{ ( l )
} ) $, and $\mathbf{B}_{l} ( \cdot ) = (
\mathbf{B}%
_{l1} ( \cdot )$, $\ldots,\mathbf
{B}_{lJ_{n}} ( \cdot )
) ^{\T}$. The empirical norm of the spline function $s_{l}$ is
\[
\Vert\bolds{\gamma}_{l}\Vert_{\mathbf{K}_{l}}= \Biggl\{
\frac
{1}{n}%
\sum_{i=1}^{n}
\frac{1}{T}\sum_{t=1}^{T}s_{l}^{2}
\bigl( \mathbf {X}%
_{it}^{ ( l ) } \bigr) \Biggr\}
^{1/2}=\Vert s_{l}\Vert_{n}.
\]
The advantage of choosing the penalization using $\Vert
s_{l}\Vert_{n}$ is that it no longer relies on a particular choice of
spline bases. This type of penalization ensures that the coefficients within
the same nonparametric component are treated as an entire group in model
selection and, therefore, it achieves the same effect as the group-wise model
selection approach \cite{yuan:lin:06}.

The penalty function $p_{\lambda_{n}} ( \cdot ) $ can be the
$%
L_{1}$-penalty with $p_{\lambda_{n}} ( \llvert \cdot\rrvert
) =\lambda_{n}\llvert \cdot\rrvert $ which provides a
LASSO estimator, or the $L_{2}$ penalty $p_{\lambda_{n}} (
\llvert \cdot\rrvert  ) =\lambda_{n}\llvert \cdot
\rrvert ^{2}$ which produces a ridge-type estimator. However, we do
not apply the $L_0$ penalty here as it is highly computationally
intensive and unstable. The smoothly clipped absolute deviation (SCAD)
\cite{fan:li:01} penalty is considered here,
where the derivative is defined as
\[
p_{\lambda_{n}}^{\prime} ( \theta ) =\lambda_{n} \biggl\{ I (
\theta\leq\lambda_{n} ) +\frac{ ( a\lambda_{n}-\theta
)
_{+}}{ ( a-1 ) \lambda_{n}}I ( \theta>\lambda
_{n} ) \biggr\},
\]
here the constant $a$ is chosen to be $3.7$ as in \cite{fan:li:01},
and $\lambda_{n}>0$ is a tuning parameter, whose selection is
described in Section~\ref{subsec:tunning}. The SCAD penalty has several
advantages such as unbiasedness, sparsity and continuity.

The penalized estimator $\widehat{\bolds{\gamma}}^{\pqif}$ is
obtained by
minimizing the penalized objective function in (\ref{DEF:pqif}). Then
for any $\mathbf{x}\in{}[0,1]^{d_{x}}$, the estimator of the
nonparametric functions in (\ref{mgaplm}) is calculated by
\[
\widehat{\alpha}_{l}^{\pqif} \bigl( x^{ ( l ) } \bigr) =
\sum_{j=1}^{J_{n}}\widehat{
\gamma}_{lj}^{\pqif}B_{lj} \bigl( x^{ (
l ) }
\bigr),\qquad l=1,\ldots,d_{x}.
\]

We establish the asymptotic properties of the penalized parametric and
nonparametric
components estimators for the marginal GAPLM in the following theorems.
We assume that in the true model only the first $s_z$ $(0\leq
s_z\leq d_z)$ linear components and the first $s_x$ $(0\leq s_x\leq d_x)$
nonlinear components are nonzero, and the remaining components are all
zeros. Let $\alpha_{0} ( \mathbf{x}_{it} ) =\sum_{l=1}^{d_{x}}\alpha
_{0,l}(x_{it}^{(l)})=\sum_{l=1}^{s_{x}}\alpha_{0,l}( x_{it}^{(l)})
+\sum_{l=s_{x}+1}^{d_{x}}\alpha_{0,l}( x_{it}^{(l)})$, with $\alpha
_{0,l}=0 $ almost surely for $l=s_{x}+1,\ldots,d_{x}$, where $s_{x}$
is the
number of nonzero nonlinear components. Similarly, let $s_{z}$ be the
number of nonzero components of $\bolds{\beta}_{0}$. Let
$\bolds{\beta}%
_{0}=(\beta_{0,1},\ldots,\beta_{0,d_{z}})^{\T}=(\bolds{\beta
}_{\mathcal{%
S}0}^{\T},\bolds{\beta}_{\mathcal{N}0}^{\T})^{\T}$, where
$\bolds{\beta}%
_{\mathcal{S}0}$ consists of all $s_{z}$ nonzero components of
$\bolds{%
\beta}_{0}$, and $\bolds{\beta}_{\mathcal{N}0}=\mathbf{0}$
without loss
of generality. In a similar fashion to $\bolds{\beta}_{0}$,
denote $%
\widehat{\bolds{\beta}}^{\pqif}= \{(\widehat{\bolds{\beta}}_{%
\mathcal{S}}^{\pqif})^{\T},(\widehat{\bolds{\beta
}}_{\mathcal
{N}}^{\pqif%
})^{\T} \} ^{\T}$.

We first derive the convergence rate of the penalized QIF estimators
$\widehat{\bolds{\beta}}^{\pqif}$ and $\{\widehat{\alpha
}_{l}^{\pqif}\}_{l=1}^{d_{x}}$. In particular, if $d_x$ is finite, we
show that this convergence rate is the same as the rate of convergence
for the unpenalized estimators $\widehat{\bolds{\beta}}^{\qif}$
and $\{\widehat{\alpha}_{l}^{\qif}\}_{l=1}^{d_{x}}$\vspace*{2pt} in {Theorem~\ref
{THM:SCAD-QIF}}. Furthermore, we prove that the penalized estimators
$\widehat{\bolds{\beta}}^{\pqif}$, $\{\widehat{\alpha
}^{\pqif}\}
_{l=1}^{d_{x}}$ possess the sparsity property as in {Theorem~\ref
{THM:sparsity}}. That is, $\widehat{\alpha}_{l}^{\pqif}=0$ almost
surely for $l=s_{x}+1,\ldots,d_{x}$, and $\widehat{\bolds{\beta
}}_{\mathcal{N}} ^{\qif}=0$. The sparsity property
implies that the model selection procedure is consistent, that is, the
selected model converges to the corrected model with probability
tending to one. We define
%
\begin{eqnarray}\label{DEF:anbn}
a_{n}&=&\max_{1\leq l\leq d_{z}}\bigl\{\bigl|p_{\lambda_{2,l}}^{\prime}
\bigl(|\beta _{0,l}|\bigr)\bigr|, \beta_{0,l}\neq0\bigr\},
\nonumber
\\[-8pt]
\\[-8pt]
\nonumber
b_{n}&=&
\max_{1\leq l\leq d_{z}} \bigl\{\bigl|p_{\lambda_{2,l}}^{\prime\prime}\bigl(|
\beta_{0,l}|\bigr)\bigr|,\beta_{l0}\neq 0\bigr\}.
\end{eqnarray}

\begin{theorem}
\label{THM:SCAD-QIF} Under conditions \textup{(C1)--(C9)} and \textup{(P2)} in Appendix %
\ref{subsec:assu}, if $d_x/ \log(n)$ $\rightarrow0$,
$n^{-1/4}d_z\rightarrow0$, $J_{n}\asymp n^{b}$, for some $1/(4r)\leq b
< 1/4$ with smoothing parameter $r>1$ defined in
condition \textup{(C1)}, and the tuning parameters $\lambda_{jn}\rightarrow0$,
$j=1,2$, $n\rightarrow\infty$, then there exists a local solution
$\widehat{\bolds{\beta}%
}^{\pqif}$ in (\ref{DEF:pqif}) such that its rate of convergence is
$O_{P}\{%
{ \rho_n^{-3}d_n^{1/2}(n^{-1/2}+a_{n})}\}$, and there exists
a local minimizer of (\ref{DEF:pqif}) such that
\[
\frac{1}{n}\sum_{l=1}^{d_x}\sum
_{i=1}^{n}\sum
_{t=1}^{T} \bigl\{ \widehat{ 
\alpha}_{l}^{\pqif}\bigl( x_{it}^{(l)}\bigr)
-\alpha_{0,l}\bigl(x_{it}^{(l)}\bigr) \bigr
\}^{2} =O_{P} \bigl\{{ \rho_n^{-6}}d_n
\bigl(n^{-1/2}+a_{n}\bigr)^2%
\bigr\}.
\]
\end{theorem}

\begin{remark}
If the number of nonparametric functions, $d_x$, is finite, then $\rho
_n$ is a fixed constant. Further if $J_n\asymp n^{1/(2r+1)}$ and
$a_n=O(n^{-1/2})$, we obtain the optimal nonparametric convergence rate
$n^{-2r/ ( 2r+1 )}$ as in \cite{xue:qu:zhou:10}. For the
parametric terms, if $d_x$ is finite, and $n^{1/(4r)}\ll J_n\ll d_z\ll
n^{1/4}$, then we obtain the same parametric convergence rate as in
\cite{cho:qu:13}.
\end{remark}

\begin{theorem}
\label{THM:sparsity} Assume that conditions \textup{(C1)--(C9)}, \textup{(P1)--(P2)} in
Appendix \ref{subsec:assu} hold. If
{ $d_x/ \log(n)\rightarrow0$, $n^{-1/5}d_z\rightarrow0$,
$J_{n}\asymp n^{b}$, for some $1/(4r)\leq b < 1/5$ with smoothing
parameter $r>1$ defined in
condition \textup{(C1)}}, and the tuning parameters $\lambda_{jn}\rightarrow
0$, and
${ \rho_n^{-1}d_n^{-1/2}n^{1/2}\lambda
_{jn}}\rightarrow\infty$, $j=1,2$, $n\rightarrow\infty$, then with
probability approaching $1$, $\widehat{\alpha}_{l}=0$ almost surely
for $
l=s_{x}+1,\ldots,d_{x}$, and the estimator {$\widehat{\bolds{\beta}}^{%
\pqif}$} has the sparsity property, that is, $P(\widehat{\bolds{\beta}}_{%
\mathcal{N}}^{\pqif}=0)\rightarrow1$ as $n\rightarrow\infty$.
\end{theorem}

Theorem~\ref{THM:sparsity} indicates that the proposed selection method
possesses model selection consistency. Theorems \ref{THM:SCAD-QIF} and
\ref{THM:sparsity} provide similar results for the nonparametric
components as those for the penalized generalized additive models in
\cite{xue:qu:zhou:10} when $d_x$ is finite. However, the theoretical
proof is very different from the penalized generalized additive model
approach and is much more challenging, due to the involvement of both
parametric and nonparametric components, where two sets of covariates
could be dependent, and the dimensions of linear and nonlinear terms
increase along with the sample size.

We also investigate the asymptotic distribution of the estimators for\break  the
parametric term. Define a vector $\bolds{\kappa}_{%
\mathcal{S}} = \{p_{\lambda_{2,n}}^{\prime}(|\beta_{0,1}|)\sgn
(\beta
_{0,1}),\ldots,\break  p_{\lambda_{2,n}}^{\prime}(|\beta_{0,s_{z}}|)
\times \sgn(\beta
_{0,s_{z}})\}^{\T}$ and a diagonal matrix $\bolds{\Lambda
}_{\mathcal{S}}=%
\diag\{p_{\lambda_{2,n}}^{\prime\prime}(|\beta_{0,1}|),\ldots,\break 
p_{\lambda_{2n}}^{\prime\prime}(|\beta_{0,s_{z}}|)\}$.\vadjust{\goodbreak} In a similar
fashion to $%
\bolds{\beta}$, we write the collections of all components,
$\mathbf{X}_{i}=(\mathbf{X}_{\mathcal{S}i}^{\T},\mathbf
{X}_{\mathcal{N}i}^{\T%
})^{\T}$, $\mathbf{Z}_{i}=(\mathbf{Z}_{\mathcal{S}i}^{\T
},\mathbf{Z}_{%
\mathcal{N}i}^{\T})^{\T}$, $\widehat{\mathbf{Z}}_{i}=(\widehat
{\mathbf{Z}}_{%
\mathcal{S}i}^{\T},\widehat{\mathbf{Z}}_{\mathcal{N}i}^{\T
})^{\T}$.
Further denote $%
\widehat{\mathbf{J}}_\DZcalS^{\T}= \{ (\widehat
{\mathbf
{J}}_\DZcalS%
^{(1)})^{\T},\ldots, (\widehat{\mathbf{J}}_\DZcalS^{(K)})^{\T
} \}
^{\T}$,
where $\widehat{\mathbf{J}}_\DZcalS^{(k)}=\frac{1}{n}\sum_{i=1}^{n}%
\mathbf{D}_{i}^{\T}\bolds{\Gamma}_{0,i}^{ ( k )
}\widehat{\mathbf{%
Z}}_{\mathcal{S}i}$. Next, let $\widehat{\bolds{\Psi
}}_{\mathcal
{S},n}=\widehat{\mathbf{J}}_{\DZcalS}^{\T%
} ( \mathbf{C}_{n}^{0} )^{-1}\widehat{\mathbf
{J}}_\DZcalS$, $
\widehat{\bolds{\Omega}}_{\mathcal{S},n}=\frac{1}{n}\sum_{i=1}^{n} \{
\widehat{\mathbf{J}}_\DZcalS^{\T} ( \mathbf
{C}_{n}^{0}
) ^{-1}%
\mathbf{W}_{i} \} ^{\otimes2}$ with $\mathbf{W}_{i}$ in
(\ref{DEF:Wi}).

\begin{theorem}
\label{THM:oracle} Assume conditions \textup{(C1)--(C9)}, \textup{(P1)--(P2)}
in Appendix
\ref%
{subsec:assu} hold.
{ If $d_x/ \log(n)\rightarrow0$, $n^{-1/5}d_z\rightarrow0$,
$J_{n}\asymp n^{b}$, for some $1/(2r+1)\leq b < 1/5$ with smoothing
parameter $r>2$ in
condition \textup{(C1)}}, and the tuning parameters $\lambda_{jn}\rightarrow
0$, $%
{ \rho_n^{-1}d_n^{-1/2}n^{1/2}\lambda
_{jn}}\rightarrow+\infty$, $j=1,2$, as $n\rightarrow\infty$, then
\[
\sqrt{n}\mathbf{A}_{n}\bolds{\Sigma}_{\mathcal
{S},n}^{-1/2}(
\widehat{%
\bolds{\Psi}}_{\mathcal{S},n}+\bolds{
\Sigma}_{\mathcal
{S}})\bigl\{\bigl(\widehat{%
\bolds{
\beta}}_{\mathcal{S}}^{\pqif}-\bolds{\beta }_{\mathcal
{S}0}
\bigr)+(%
\widehat{\bolds{\Psi}}_{\mathcal{S},n}+\bolds{
\Lambda }_{\mathcal{S}%
})^{-1}\bolds{\kappa}_{\mathcal{S}}\bigr
\} \mathop{\rightarrow}^{D} N(0,\bolds{\Sigma}%
_{A}),
\]
where $\mathbf{A}_n$ is any $q \times d_z$ matrix with a finite
$q$ such
that $\bolds{\Sigma}_{A}=\lim_{n\rightarrow\infty}
\mathbf{A}
_{n}^{\otimes2}$, and $\bolds{\Sigma}_{\mathcal{S},n}=\widehat
{\bolds{%
\Psi}}_{\mathcal{S},n}^{-1}{\widehat{\bolds{\Omega
}}}_{\mathcal
{S},n}%
\widehat{\bolds{\Psi}}_{\mathcal{S},n}^{-1}$.
\end{theorem}

\subsection{An algorithm}
\label{subsec:LQA}

To minimize the PQIF in (\ref{DEF:pqif}), we develop an algorithm based
on the local quadratic approximation \cite{fan:li:01}. To obtain an
initial estimator $ ( \bolds{\beta}^{0},\bolds{\gamma
}^{0} ) $
which is sufficiently close to the true minimizer of (\ref{DEF:pqif}), we
could choose the unpenalized QIF estimator $\widehat{%
\bolds{\theta}}^{\qif}=\{(\widehat{\bolds{\beta}}^{\qif
})^{\T},(\widehat{%
\bolds{\gamma}}^{\qif})^{\T}\}^{\T}$ as the initial value. Let
$\bolds{%
\beta}^{k}= ( \beta_{1}^{k},\ldots,\beta_{d_{z}}^{k} )
^{\T}$
and $\bolds{\gamma}^{k}= ( \bolds{\gamma}_{1}^{k\T
},\ldots,\bolds{%
\gamma}_{d_{x}}^{k\T} ) ^{\T}$ be the values at the $k$th iteration.
If $\beta_{l}^{k}$ (or $\bolds{\gamma}_{l^{\prime}}^{k}$) is
close to
zero, such that $|\beta_{l}^{k}|\leq\epsilon$ (or $\Vert\bolds{\gamma}%
_{l^{\prime}}^{k}\Vert_{\mathbf{K}_{l^{\prime}}}\leq\epsilon$)
with some
small threshold value $\epsilon$, then $\beta_{l}^{k+1}$ (or
$\bolds{%
\gamma}_{l^{\prime}}^{k+1}$ ) is set to $\mathbf{0}$. We consider
$\epsilon=10^{-6}$ in our numerical examples.\vspace*{1pt}

Suppose $\beta_{l}^{k+1}=0$, for $l=b_{k}+1,\ldots,d_{z}$, and
$\bolds{%
\gamma}_{l}^{k+1}=0$, for $l=b_{k}^{\prime}+1,\ldots,d_{x}$, and
$\bolds{%
\theta}^{k+1}  = (\beta_{1}^{k+1},\ldots, \beta
_{b_{k}}^{k+1},\beta
_{b_{k}+1}^{k+1},\ldots,\beta_{dz}^{k+1},( \bolds{\gamma
}_{1}^{k+1}) ^{\T},\ldots, ( \bolds{\gamma}_{b_{k}^{\prime}}^{k+1}) ^{\T},
( \bolds{\gamma
}_{b_{k}^{\prime}+1}^{k+1}) ^{\T},\break \ldots,  ( \bolds{\gamma
}_{d_{x}}^{k+1})^{\T}) ^{\T}=\{( \bolds{\beta}_{\mathcal{S}}^{k+1}) ^{\T}, (
\bolds{\beta}_{\mathcal{N}}^{k+1})^{\T}, (\bolds{\gamma}_{\mathcal{S}}^{k+1})
^{\T}, ( \bolds{\gamma}_{\mathcal{N}}^{k+1}) ^{\T}\}^{\T}$, in
which $\bolds{\beta}_{\mathcal{N}}^{k+1}=\mathbf{0}$,
$\bolds{\gamma}_{\mathcal{N}}^{k+1}=\mathbf{0}$. Let
$\bolds{\theta}%
= ( \bolds{\beta}_{\mathcal{S}}^{\T},\bolds{\beta
}_{\mathcal{N}}^{\T%
},\bolds{\gamma}_{\mathcal{S}}^{\T},\bolds{\gamma
}_{\mathcal
{N}}^{\T%
} ) ^{\T}$ be the partition of any $\bolds{\theta}$.

The local quadratic\vspace*{1pt} approximation is implemented for obtaining the nonzero
components $\bolds{\theta}_{\mathcal{S}}^{k+1}=\{(\bolds{\beta}_{%
\mathcal{S}}^{k+1})^{\T}, ( \bolds{\gamma}_{\mathcal{S}}^{k+1})
^{\T}\}^{\T%
}$. Specifically, for $|\beta_{^{l}}^{k}|>\epsilon$, the penalty
for the parametric term is approximated by
\begin{eqnarray*}
p_{\lambda_{n}} \bigl( |\beta_{l}| \bigr) &\approx&p_{\lambda
_{n}} \bigl( \bigl|
\beta_{l}^{k}\bigl| \bigr) +p_{\lambda_{n}}^{\prime}
\bigl(\bigl |\beta _{l}^{k}\bigr| \bigr) \bigl( \llvert
\beta_{l}\rrvert -\bigl \vert \beta _{l}^{k} \bigr\vert
\bigr)
\\
&\approx&p_{\lambda_{n}} \bigl( \bigl|\beta_{l}^{k}\bigr| \bigr) +
\tfrac{1}{2} 
p_{\lambda_{n}}^{\prime} \bigl( \bigl|
\beta_{l}^{k}\bigr| \bigr) \bigl|\beta _{l}^{k}\bigr|^{-1}
\bigl\{ \beta_{l}^{2}-\bigl( \beta_{l}^{k}
\bigr) ^{2} \bigr\}.
\end{eqnarray*}
For $\Vert\bolds{\gamma}_{l^{\prime}}^{k}\Vert_{\mathbf
{K}%
_{l^{\prime}}}>\epsilon$, the penalty function for the nonparametric part
is approximated by
\begin{eqnarray*}
&& p_{\lambda_{n}} \bigl(\Vert\bolds{\gamma }_{l^{\prime
}}
\Vert_{\mathbf{K}%
_{l^{\prime}}} \bigr)
\\
&&\qquad  \approx p_{\lambda_{n}} \bigl(\bigl\Vert\bolds{\gamma}%
_{l^{\prime}}^{k}\bigr\Vert_{\mathbf{K}_{l^{\prime}}} \bigr) +p_{\lambda
_{n}}^{\prime}
\bigl(\bigl \Vert\bolds{\gamma}_{l^{\prime
}}^{k}\bigr\Vert
_{\mathbf{%
K}_{l^{\prime}}} \bigr) \bigl\Vert\bolds{\gamma}_{l^{\prime
}}^{k}
\bigr\Vert_{%
\mathbf{K}_{l^{\prime}}}^{-1}\bolds{\gamma}_{l^{\prime
}}^{kT}
\mathbf{K}%
_{l^{\prime}} \bigl( \bolds{
\gamma}_{l^{\prime}}-\bolds{\gamma}%
_{l^{\prime}}^{k}
\bigr)
\\
&&\qquad \approx p_{\lambda_{n}} \bigl(\bigl\Vert\bolds{\gamma }_{l^{\prime
}}^{k}
\bigr\Vert _{\mathbf{K}_{l^{\prime}}} \bigr) +\tfrac{1}{2}p_{\lambda
_{n}}^{\prime
}
\bigl(\bigl\Vert\bolds{\gamma}_{l^{\prime}}^{k}\bigr\Vert
_{\mathbf{K}
_{l^{\prime}}} \bigr) \bigl\Vert\bolds{\gamma}_{l^{\prime
}}^{k}
\bigr\Vert _{\mathbf{%
K}_{l^{\prime}}}^{-1} \bigl( \bolds{\gamma}_{l^{\prime}}^{\T
}
\mathbf{K}%
_{l^{\prime}}\bolds{\gamma}_{l^{\prime}}-
\bolds{\gamma }_{l^{\prime}}^{kT}%
\mathbf{K}_{l^{\prime}}\bolds{\gamma}_{l^{\prime
}}^{k}
\bigr),
\end{eqnarray*}
where $p_{\lambda_{n}}^{\prime}$ is the first-order derivative of $%
p_{\lambda_{n}}$.\vadjust{\goodbreak}

This leads to the local approximation of the objective function in
(\ref
{DEF:pqif}) by a quadratic function:
\begin{eqnarray*}
&&Q_{n}\bigl( \bolds{\theta}^{k}\bigr) +
\dot{Q}_{n}\bigl( \bolds{\theta }^{k}
\bigr)^{\T} \pmatrix{ \bolds{
\beta}_{\mathcal{S}}-\bolds{\beta}_{\mathcal{S}}^{k}
\vspace*{2pt}\cr
\bolds{\gamma}_{\mathcal{S}}-\bolds{\gamma}_{\mathcal
{S}}^{k}}
 +\frac{1}{2}\pmatrix{ \bolds{\beta}_{\mathcal{S}}-\bolds{
\beta}_{\mathcal{S}}^{k}
\vspace*{2pt}\cr
\bolds{\gamma}_{\mathcal{S}}-\bolds{\gamma}_{\mathcal
{S}}^{k}} ^{\T}\ddot{Q}_{n}\bigl(
\bolds{\theta}^{k}\bigr) \pmatrix{
\bolds{\beta}_{\mathcal{S}}-\bolds{\beta}_{\mathcal{S}}^{k}
\vspace*{2pt}\cr
\bolds{\gamma}_{\mathcal{S}}-\bolds{\gamma}_{\mathcal
{S}}^{k}}
\\
&&\qquad{}+\frac{1}{2}n\pmatrix{
\bolds{
\beta}_{\mathcal{S}}
\vspace*{2pt}\cr
\bolds{\gamma}_{\mathcal{S}}}
^{\T}\Lambda\bigl(\bolds{\theta}^{k}\bigr) \pmatrix{
\bolds{\beta}_{\mathcal{S}}
\vspace*{2pt}\cr
\bolds{\gamma}_{\mathcal{S}}%
},
\end{eqnarray*}
where $\dot{Q}_{n} ( \bolds{\theta}^{k} ) =\frac
{\partial
Q_{n} ( \bolds{\theta}^{k} ) }{\partial\bolds{\theta}_{%
\mathcal{S}}}$, $\ddot{Q}_{n} ( \bolds{\beta}^{k} )
=\frac{%
\partial^{2}Q_{n} ( \bolds{\theta}^{k} ) }{\partial
\bolds{%
\theta}_{\mathcal{S}}\,\partial\bolds{\theta}_{\mathcal
{S}}^{\T
}}$ with $%
\bolds{\theta}_{\mathcal{S}}= ( \bolds{\beta
}_{\mathcal
{S}}^{\T},%
\bolds{\gamma}_{\mathcal{S}}^{\T} ) ^{\T}$, and
%
\begin{eqnarray}\label
{equ_Lambda}\qquad
&&\Lambda \bigl( \bolds{\theta}^{k} \bigr)=\diag
\bigl\{\bigl|\beta_{1}^{k}\bigr|^{-1} p_{\lambda_{n}}^{\prime}
\bigl(\bigl|\beta_{1}^{k}\bigr|\bigr),\ldots, \bigl|\beta_{b_{k}}^{k}\bigr|^{-1}
p_{\lambda_{n}}^{\prime}\bigl(\bigl|\beta_{b_{k}}^{k}\bigr|
\bigr),
\nonumber
\\[-4pt]
\\[-12pt]
\nonumber
&&\hspace*{43pt}\qquad \bigl\Vert\bolds{\gamma}_{1}^{k}
\bigr\Vert_{\mathbf{K}1}^{-1} p_{\lambda_{n}}^{\prime}\bigl(\bigl\Vert
\bolds{\gamma}_{1}^{k}\bigr\Vert_{%
\mathbf{K}_{1}}\bigr)
\mathbf{K}_{1},\ldots, \bigl\Vert\bolds{%
\gamma}_{b_{k}^{\prime}}^{k}\bigr\Vert_{\mathbf{K}_{b_{k}^{\prime}}}^{-1}
p_{\lambda_{n}}^{\prime}\bigl(\bigl\Vert\bolds{\gamma}%
_{b_{k}^{\prime}}^{k}\bigr\Vert_{\mathbf{K}_{b_{k}^{\prime
}}}\bigr)\mathbf{K%
}_{b_{k}^{\prime}} \bigr\}.
\end{eqnarray}
We minimize the above quadratic function to get $\bolds{\theta
}_{\mathcal{S}}^{k+1}$. The corresponding Newton--Raphson algorithm provides
\[
\bolds{\theta}_{\mathcal{S}}^{k+1}=\bolds{\theta
}_{\mathcal
{S}%
}^{k}- \bigl\{ \ddot{Q}_{n} \bigl(
\bolds{\theta}^{k} \bigr) +n\Lambda \bigl( \bolds{
\theta}^{k} \bigr) \bigr\} ^{-1} \bigl\{ \dot
{Q}_{n} \bigl( \bolds{\theta}^{k} \bigr) +n\Lambda
\bigl( \bolds{\theta }^{k} \bigr) \bolds{
\theta}_{\mathcal{S}}^{k} \bigr\}.
\]
The above iteration process is repeated until convergence is reached,
where the convergence criterion is based on $\llVert \bolds{\theta}%
^{k+1}-\bolds{\theta}^{k}\rrVert \leq10^{-6}$. The
proposed algorithm
is quite stable and converges quickly. However, in general, the
computational time increases as the
dimension of covariates increases.

\subsection{Tuning parameter and knots selection}
\label{subsec:tunning}

Tuning parameter and knots selections play important roles in the
performance of model selection. The spline approximation for the
nonparametric components requires an appropriate selection of the knot
sequences $ \{ \upsilon_{l} \} _{l=1}^{d_x}$ in Section~\ref{subsec:splineQIF}. For the penalized QIF method in Section~\ref{subsec:pQIF}, in addition to knots selection, we also need to
address how to choose tuning parameters $\lambda_{1,n}$ and $\lambda
_{2,n}$ in the SCAD penalty function. To reduce computational
complexity, we consider $ \lambda_{1,n}=\lambda_{2,n}=\lambda_{n}$
and select only $\lambda_{n}$. This is justified by Theorems \ref
{THM:SCAD-QIF}, \ref{THM:sparsity} and \ref{THM:oracle} in Section~\ref{sec:pQIF}.

Although selecting the number and position of spline knots is important
in curve smoothing, in our simulation study we found that knot
selection seems to be less critical for the estimation of the
parametric coefficients and model selection than for the estimation of
the nonparametric components. For convenience, we choose equally spaced
knots and the number of interior knots is selected as the integer part
of $N_n=n^{1/(2p+3)}$, where $n$ is the sample size and $p$ is the
order of the polynomial spline. This approach is also adopted in \cite
{huang:wu:zhou:04,xue:qu:zhou:10} and \cite{xue:qu:12}. Furthermore,
we use the same knot sequences selected in the unpenalized procedure
for the penalized QIF estimation. Therefore, we only need to determine
the tuning parameter for the penalization part. For any given tuning
parameter $\lambda_{n}$, the estimator minimizing~(\ref{DEF:pqif}) is
denoted as $\widehat{\bolds{\theta}}_{\lambda_{n}}$. We propose
to use the extended Bayesian Information Criterion~(EBIC) to select the
optimal tuning parameters based on \cite{chen:chen:08} and \cite
{huang:horowitz:wei:10}. Because the QIF $Q_{n}$ is analog to minus
twice the log-likelihood function \cite{qu:lindsay:li:00}, we define
the EBIC in the PQIF procedure as
%
\begin{eqnarray}\label{EBIC1}
\operatorname{EBIC} ( \lambda_{n} ) & = & Q_{n}(\widehat {
\bolds{\theta}}%
_{\lambda_{n}})+\log ( n ) \hat{d}_{z}(
\lambda_{n})+\log \bigl( \nu_{z}(\lambda_{n})
\bigr)
\nonumber
\\[-8pt]
\\[-8pt]
\nonumber
&&{}+\log ( n ) N_{n}\hat{d}_{x}(\lambda
_{n})+N_{n}\log \bigl( \nu_{x}(
\lambda_{n}) \bigr),
\nonumber
\end{eqnarray}
where $\hat{d}_{z}(\lambda_{n})$ and $\hat{d}_{x}(\lambda_{n})$ are the
nonzero parametric and nonparametric terms in $\widehat{\bolds{\theta}}%
_{\lambda_{n}}$, respectively, and $\nu_{z}(\lambda_{n})=
{d_{z}\choose
\hat{d}_{z}(\lambda_{n})}$, which is a combination operator and
represents the number of choices for selecting $\hat{d}_{z}(\lambda
_{n})$ terms out of $d_{z}$ parametric terms. Similarly, define $\nu
_{x}(\lambda_{n})={d_{x}\choose\hat{d}_{x}(\lambda_{n})}$. %
See \cite{chen:chen:08} for details. However, when the full likelihood
is available, it is more accurate to use minus twice the log-likelihood
function instead of $Q_{n}$ as the first term in (\ref{EBIC1}). That is,
\begin{eqnarray*}\label{EBIC}
\operatorname{EBIC} ( \lambda_{n} ) & = & -2\log L(\widehat {\bolds{
\theta}}%
_{\lambda_{n}})+\log ( n ) \hat{d}_{z}(
\lambda_{n})+\log \bigl( \nu_{z}(\lambda_{n})
\bigr)
\\
&&{}+\log ( n ) N_{n}\hat{d}_{x}(\lambda
_{n})+N_{n}\log \bigl( \nu_{x}(
\lambda_{n}) \bigr),
\end{eqnarray*}
where $L(\cdot)$ is the full likelihood function. As indicated in
\cite
{wang:qu:09}, the one using the full likelihood, if it is available,
has better finite sample performance when the sample size is small. The
optimal $\lambda_{n}$ is chosen such that the EBIC value reaches the
minimum, or equivalently, $\widehat{\lambda}_{n}=\argmin_{\lambda_{n}}
\operatorname{EBIC} ( \lambda_{n} ) $.

\section{Simulation studies}\label{sec:simu}

In this section, we assess the numerical performance of the proposed
methods through simulation studies. To assess estimation accuracy and
efficiency, define the
model error (ME) as
\[
\frac{1}{n^{*}T}\sum_{i=1}^{n^{*}}\sum
_{t=1}^{T} \Biggl\{g^{-1} \Biggl(
\sum_{l=1}^{d_x}\widehat{
\alpha}_{l} \bigl(x_{it}^{(l)} \bigr)+
\mathbf{z}%
_{it}^{\T}\widehat{\bolds{
\beta}} \Biggr) -g^{-1} \Biggl(\sum_{l=1}^{d_x}
\alpha_{l} \bigl(x_{it}^{(l)} \bigr)+
\mathbf{z}%
_{it}^{\T}\bolds{\beta} \Biggr)
\Biggr\}^{2},
\]
where $ (\mathbf{x}_{it}, \mathbf{z}_{it}
)_{i=1,t=1}^{n^{*},T}$ are independently generated test data and follow
the same distribution as the training data. In our simulations, we take
$n^{*}=1000$. Furthermore, $g^{-1}$ is the identity link function for
continuous outcomes and the logit link function for binary outcomes.
The model error measures the prediction performance of different
methods. Denote the index sets of the selected and true models by $\hat
{\mathcal{S}}$ and $\mathcal{S}_{0}$, respectively. If $\hat
{\mathcal
{S}}=\mathcal{S}_{0}$, then $\hat{\mathcal{S}}$ is a correct selection;
if $\mathcal{S}_{0}\subset\hat{\mathcal{S}}$ and $\mathcal
{S}_{0}\neq
\hat{\mathcal{S}}$, then we call $\hat{\mathcal{S}}$ over selection;
otherwise, if $\mathcal{S}_{0}%
\not\subset\hat{\mathcal{S}}$, then $\hat{\mathcal{S}}$ under
selection. The number of replications is $500$ in the following
simulation studies.

\subsection{Example~1: Continuous response}

The continuous responses $ \{Y_{it} \}$ are generated from
%
\begin{equation}
Y_{it}=\sum_{l=1}^{d_x}
\alpha_{l} \bigl(X_{it}^{(l)} \bigr)+\mathbf
{Z}_{it}^{%
\T}\bolds{\beta}+\varepsilon_{it},\qquad
 i=1,\ldots,n, t=1,\ldots,5, \label{Simu}
\end{equation}
where $n=100,200$, or $500$, and $d_x=d_z=2n^{1/4}$ which is rounded to the
nearest integer and takes values of $6,8$ and $10$, respectively, for $%
n=100,200$ and $500$. We take $\alpha_{1} ( x )
=\sin (2\pi x )$, $\alpha_{2} ( x ) =8x(1-x)-4/3$,
and $%
\alpha_{l} ( x ) =0$ for $l=3,\ldots,d_x$, and $\beta
_{1}=1,\beta_{2}=2$, and $\beta_{l}=0$ for $l=3,\ldots,d_z$. Therefore,
only the
first two variables in $\mathbf{X}_{it}$ and $\mathbf{Z}_{it}$
are relevant and the rest are null variables. The covariates
$\mathbf{X}_{it}=(
X_{it}^{(1)},\ldots,X_{it}^{(d_x)}) ^{\T}$ are generated by $%
X_{it}^{(l)}= (2W_{it}^{(l)}+U_{it} )/3$, where $\mathbf{W}
_{it}= (W_{it}^{(1)},\ldots,W_{it}^{(d_x)} )$
and $U_{it}$ are independently generated from
$\operatorname{Uniform} ( [0,1]^{d_x} ) $ and $\operatorname{Uniform} ( [0,1] )$,
respectively.
Therefore, the covariates $\mathbf{X}_{it}$ have an exchangeable
correlation structure. In addition, $\mathbf{Z}_{it}=(
Z_{it}^{(1)},\ldots
,Z_{it}^{(d_z)}) ^{\T}$ are generated with $Z_{it}^{(1)}=1$ and $
(Z_{it}^{(2)},\ldots,Z_{it}^{(d_z)} )$ being
generated from a zero mean multivariate normal distribution with a
marginal variance of
$1$ and an AR-1 correlation with parameter $\rho=0.7$. The errors
$\bolds{%
\varepsilon}_{i}= ( \varepsilon_{i1},\ldots,\varepsilon
_{i5} )
^{\T}$ follows a zero mean multivariate normal with a
marginal variance of $\sigma^{2}=1.5$ and an exchangeable correlation
with correlation $\rho=0.7$.

The proposed penalized QIF method with the SCAD penalty is considered.
In spline approximation, we use both the linear splines and cubic
splines. Furthermore, we consider basis matrices from three different
working correlation structures: exchangeable (EC), AR-1 and independent
(IND), and compare their estimation efficiencies to illustrate the
effect on efficiency gain of incorporating within-cluster correlation.

\begin{table}
\caption{Example~1: The simulation results using the SCAD penalty with
exchangeable (EC), AR-1 or independent (IND) working correlation and
linear or cubic splines. The columns of C, O and U provide the
percentage of correct selection, over selection and under selection,
and MME provides the averaged model errors from 500 replications}
\label{Guassian_low}
\begin{tabular*}{\textwidth}{@{\extracolsep{\fill}}lccccccccc@{}}
\hline
& & \multicolumn{4}{c}{\textbf{Linear spline}} & \multicolumn{4}{c@{}}{\textbf{Cubic
spline}} \\[-6pt]
& & \multicolumn{4}{c}{\hrulefill} & \multicolumn{4}{c@{}}{\hrulefill} \\
\textbf{Method} & $\bolds{n}$ & \textbf{C} & \textbf{O} & \textbf{U} &
\textbf{MME} & \textbf{C} & \textbf{O} & \textbf{U} & \textbf{MME} \\
\hline
EC & $100$ & $0.936$ & $0.050$ & $0.014$ & $0.0461$ & $0.842$ & $0.006$
& $0.098$
& $0.1337$ \\
& $200$ & $0.992$ & $0.008$ & $0.000$ & $0.0258$ & $1.000$ & $0.000$ &
$0.000$ & $%
0.0175$ \\
& $500$ & $1.000$ & $0.000$ & $0.000$ & $0.0089$ & $1.000$ & $0.000$ &
$0.000$ & $%
0.0054$ \\[3pt]
AR-1 & $100$ & $0.930$ & $0.012$ & $0.058$ & $0.0660$ & $0.766$ &
$0.012$ & $0.222$
& $0.2312$ \\
& $200$ & $0.994$ & $0.00$\phantom{0} & $0.006$ & $0.0268$ & $1.000$ & $0.000$ &
$0.000$ & $%
0.0206$ \\
& $500$ & $1.000$ & $0.000$ & $0.000$ & $0.0097$ & $1.000$ & $0.000$ &
$0.000$ & $%
0.0062$ \\[3pt]
IND & $100$ & $0.836$ & $0.008$ & $0.156$ & $0.1046$ & $0.648$ &
$0.002$ & $0.350$
& $0.2707$ \\
& $200$ & $0.986$ & $0.012$ & $0.002$ & $0.0322$ & $0.994$ & $0.006$ &
$0.000$ & $%
0.0259$ \\
& $500$ & $1.000$ & $0.000$ & $0.000$ & $0.0128$ & $1.000$ & $0.000$ &
$0.000$ & $%
0.0091$ \\
\hline
\end{tabular*}
\end{table}

Table~\ref{Guassian_low} presents the variable selection and estimation
results. It summarizes the percentages of correct selection (C), over
selection (O) and under selection (U). It also gives the mean model
errors (MME) from 500 replications. Table~\ref{Guassian_low} indicates
that the probability of recovering the correct model increases to $1$
quickly and the MME decreases as the sample size increases. This
confirms the consistency theorems of variable selection and estimation
provided in Section~\ref{subsec:pQIF}. It also shows that the
procedures with a correct EC working correlation always have the
smallest MMEs and, therefore, the estimators are more efficient than
their counterparts with IND structure, which ignore within-cluster
correlation. The method with a misspecified AR-1 correlation is less
efficient than the one using the true EC structure, but is still more
efficient than assuming independent structure. Furthermore, it also
shows that the percentage of correct model-fitting using EC structure
is higher than the one using IND when the sample size is small ($n=100$).

\subsection{Example~2: Continuous response with randomly generated
correlation structure}

To assess our method in a more challenging scenario, we consider a
model similar to (\ref{Simu}), but with randomly generated correlation
structures. In particular, we assume that the dimensions of
$\mathbf
{X}$ and $\mathbf{Z}$ are $d_{x}=9$, $d_{z}=5$, respectively. As in
(\ref{Simu}), only $X_{it}^{(1)}$, $X_{it}^{(2)}$, $Z_{it}^{(1)}$ and
$Z_{it}^{(2)}$ are relevant and take the same forms as in Example~1.
Furthermore, we consider the number of clusters $n=25$ or $250$, and
cluster size $3$. The set-up of $n=25$ mimics the real data analyzed in
Section~\ref{sec:real}. The errors $ \{\bolds{\varepsilon
}_{i}= (\varepsilon_{i1},\ldots,\varepsilon_{i3} ) ^{\T
} \}_{i=1}^{25}$ independently follow a multivariate normal
distribution as in Example~1, but with a randomly generated correlation
matrix $\bolds{\Gamma}_{r}$ for each replication $r$. Let
$\bolds{\Sigma}_{1}$ be a matrix with diagonals being $1$ and all
the off-diagonals with value $0.5$, and $\bolds{\Sigma
}_{r2}=\mathbf{Q}_{r}\bolds{\Lambda}_{r}\mathbf
{Q}^{T}_{r}$ with $\mathbf{Q}_{r}$ being a randomly generated
orthogonal matrix and $\bolds{\Lambda}_{r}=\diag (\lambda
_{r1},\lambda_{r2},\lambda_{r3} )$ with $ \{\lambda
_{rj}
\}_{j=1}^{3}$ being randomly generated from $\operatorname{Uniform}[0.2,2]$. Let
$\bolds{\Sigma}_{r}=\bolds{\Sigma}_{1}+\bolds{\Sigma
}_{r2}$ and $\sigma_{r1},\ldots,\sigma_{r3}$ be the diagonal elements
of $\bolds{\Sigma}_{r}$. Let $\bolds{\Delta}_{r}=\diag
\{
\sigma^{-1/2}_{r1},\ldots,\sigma^{-1/2}_{r3} \}$. Then the
randomly generated correlation structure for the $r$th replication is
$\bolds{\Gamma}_{r}=\bolds{\Delta}_r\bolds{\Sigma
}_{r}\bolds{\Delta}_{r}$. We use this example to investigate the
performance of the QIF method in approximating the randomly generated
correlation structures.

We estimate the model using the proposed penalized QIF method with
linear spline and SCAD penalty, and assume IND, EC or AR-1 working
correlation structure. We also consider linear spline QIF estimations
of a full model (FULL) and an oracle model (ORACLE), where the full
model contains all $14$ variables while the oracle one has only the
four nonzero variables. The oracle model is not available in real data
analysis where the underlying data-generating process is unknown.

\begin{table}
\caption{Example~2: Continuous response with randomly generated
correlation. The percentage of correct selection (C), over selection
(O) and under selection (U) are provided using linear spline with the
SCAD penalty for three working correlation: exchangeable (EC), AR-1 or
independent (IND). The columns of SCAD, ORACLE and FULL report the mean
model error (MME) of the SCAD approach
and a standard linear spline estimation of the oracle model (ORACLE),
and the full model (FULL) from $500$ replications}
\label{Guassian_high}
\begin{tabular*}{\textwidth}{@{\extracolsep{\fill}}lccccccc@{}}
\hline
$\bolds{n}$ & \textbf{Method} & \textbf{C} &
\textbf{O} & \textbf{U} & \textbf{SCAD} & \textbf{ORACLE} & \textbf{FULL} \\
\hline
250&EC & $0.998$ & $0.002$ & $0.000$ & $0.0242$ & $0.0223$ & $0.0656$ \\
&AR1 & $0.990$ & $0.002$ & $0.008$ & $0.0245$ & $0.0233$ & $0.0704$ \\
&IND & $0.986$ & $0.006$ & $0.008$ & $0.0256$ & $0.0250$ & $0.0713$ \\[3pt]
\phantom{0}25&EC & $0.616$ & $0.194$ & $0.190$ & $0.5081$ & $0.3858$ & $1.6886$ \\
&AR1 & $0.566$ & $0.212$ & $0.222$ & $0.5281$ & $0.3723$ & $1.7528$ \\
&IND & $0.536$ & $0.256$ & $0.208$ & $0.5546$ & $0.3518$ & $0.7729$ \\
\hline
\end{tabular*}
\end{table}

Table~\ref{Guassian_high} summarizes variable selection performance on
correct, over and under selection percentages of the SCAD approach with
IND, EC and AR-1 working correlations and reports the mean model error
(MME) for FULL, ORACLE and SCAD when the sample size $n=25$ and $250$,
respectively. Table~\ref{Guassian_high} clearly indicates that, for a
randomly generated correlation, SCAD with an EC working correlation
still performs better than the one with IND working structure.
Furthermore,\vadjust{\goodbreak} when the sample size is large ($n=250$), the estimation
using EC always yield a smaller MME than the one with IND working
structure. It indicates that although EC is a misspecified correlation
structure, it can still improve estimation and inference performances
by incorporating some correlation in the data into the estimation. When
the sample size is small ($n=25$), the estimation using EC or AR1
working correlations of FULL and ORACLE is worse due to the extra noise
in modeling within-cluster correlation. However, the SCAD with EC or
AR1 working correlations still give smaller MMEs than SCAD with IND
correlation, due to their better performances in recovering the correct
model. Finally, Table~\ref{Guassian_high} also shows that the penalized
procedure dramatically improves estimation accuracy compared to the
un-penalized approach, with MMEs from the SCAD being very close to the
MMEs from the ORACLE model, and much smaller than the FULL model.

\begin{figure}

\includegraphics{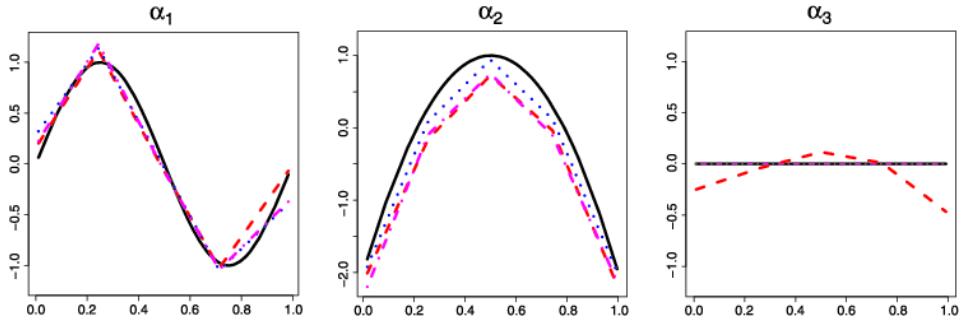}

\caption{Example~2: Plots of the first three estimated
functions from SCAD (dot--dash), Oracle (dotted) and Full (dashed)
approaches with the true functions (line). For $\alpha_3$, both SCAD
and Oracle give exactly zero estimates. The cluster size is $n=250$.}
\label{Gaussian_plot}
\end{figure}

From one selected data set, Figure~\ref{Gaussian_plot} plots the first
three estimated functional components from the SCAD, FULL and ORACLE
models using linear spline and exchangeable working correlation for
cluster size $n=250$. Note that for the third variable, both the true
and estimated functions from SCAD are zero. It shows that the proposed
SCAD method estimates unknown functions reasonably well.

\subsection{Example~3: Binary response}

A random sample of $250$ clusters is generated in each simulation run.
Within each cluster, binary responses $ \{Y_{it} \}
_{t=1}^{20}$ are generated from a marginal logit model
\[
\mathrm{logit}P (Y_{it}=1|\mathbf{X}_{it}=\mathbf
{x}_{it},\mathbf{Z}%
_{it}=
\mathbf{z}_{it} )=\sum_{l=1}^{5}
\alpha_{l} \bigl(x^{(l)}_{it}%
\bigr)+
\mathbf{z}_{it}^{\T}\bolds{\beta},
\]
where $\alpha_{1} (x ) = \cos (2\pi x )/4$,
$\alpha
_{l} (x ) = 0$, for $l=2,\ldots,5$, and $\bolds{\beta
} =
(\beta_{1},\ldots,\beta_{10} )^{\T}$ with $\beta
_{1}=1$ and
$\beta_{l}=0$ for $l=2,\ldots,10$. The covariates $\mathbf
{X}_{it}= \{X_{it}^{(l)} \}_{l=1}^{5}$ and $\mathbf
{Z}_{it}= \{Z_{it}^{(l)} \}_{l=1}^{10}$ are generated in the
same way as in Example~1. The covariates $\mathbf{X}_{it}$ have an
exchangeable correlation structure, and $\mathbf{Z}_{it}$ have an
AR-1 correlation structure with $\rho=0.7$. We use the algorithm
described in \cite{Macke:Berens:Ecker:08} to generate the correlated
binary data. It has an exchangeable correlation structure with a
correlation coefficient of $0.3$.

We conduct variable selection using the proposed penalization method
with linear spline (SCAD). We also consider estimation of the full
(FULL) and oracle (ORACLE) models using the unpenalized QIF with linear
spline. We minimize~(\ref{EQ:thetahat_qif}) and (\ref{DEF:pqif}) using
AR-1 and independent working structures, in addition to the true
exchangeable correlation structure.

\begin{table}
\caption{Example~3: Binary response. The percentages of correct
selection (C), over selection (O) and under selection (U) are provided
using linear spline with the SCAD penalty for three working
correlation: exchangeable (EC), AR-1 or independent (IND). The columns
of SCAD, ORACLE and FULL provide the mean model error ($\times10^3$)
of the SCAD approach and a standard linear spline estimation of the
oracle model (ORACLE), and the full model (FULL). The number of
replications is $500$}

\label{Binary_1}
\begin{tabular*}{\textwidth}{@{\extracolsep{\fill}}lcccccc@{}}
\hline
& \textbf{C} & \textbf{O} & \textbf{U} &
\textbf{SCAD} & \textbf{ORACLE} & \textbf{FULL} \\ \hline
EC & $0.868$ & $0.084$ & $0.048$ & $0.1102$ & $0.0547$ & $1.6820$ \\
AR1 & $0.692$ & $0.012$ & $0.292$ & $0.1264$ & $0.0586$ & $1.8489$ \\
IND & $0.684$ & $0.006$ & $0.208$ & $0.1484$ & $0.0616$ & $1.9057$ \\
\hline
\end{tabular*}
\end{table}

\begin{table}[b]
\caption{Example~3: Binary response. The sample mean and standard
deviation (in parenthesis) of $\widehat{\protect\beta}$ from the SCAD,
ORACLE and FULL model approaches}
\label{Binary_2}
\begin{tabular*}{\textwidth}{@{\extracolsep{\fill}}lccc@{}}
\hline
& \textbf{SCAD} & \textbf{ORACLE} & \textbf{FULL} \\
\hline
EC & 1.0258 (0.0461) & 1.0115 (0.0436) & 1.0945 (0.0598) \\
AR1 & 0.9969 (0.0558) & 1.0177 (0.0537) &  1.0748 (0.0738) \\
IND & 0.9932 (0.0792) & 1.0543 (0.0758) & 1.0801 (0.0893) \\
\hline
\end{tabular*}
\end{table}

Table~\ref{Binary_1} summarizes the MMEs for the SCAD, ORACLE and FULL
with three different working correlations. Table~\ref{Binary_2} also
reports the sample means and sample standard deviations (SD) of the
estimators of the nonzero regression coefficient $\widehat{\beta}_{1}$
from $500$ replications. It again shows that estimation based on
correctly specified exchangeable correlation structure is the most
efficient, having the smallest MMEs and SDs. Estimation with a
misspecified AR-1 correlation results in some efficiency loss compared
to using the true structure, but it is still much more efficient than
assuming independent structure. However, for GEE, estimation using a
misspecified AR-1 correlation structure could be less efficient than
assuming independence, since the GEE requires the estimation of the
correlation \thinspace$\rho$ for misspecified AR-1, and the estimator
of $\rho$ may not be valid.

Furthermore, similar to the previous study, MMEs calculated based on
the SCAD approach are very close to the ones from ORACLE, and much
smaller than the MMEs from the FULL model. The MMEs of the FULL model
are close to 4 times the MMEs of SCAD. This shows that the SCAD
penalization improves estimation accuracy significantly by effectively
removing the redundant variables. Table~\ref{Binary_1} also gives the
frequency of correct, over and under selection for the SCAD approach.
Overall, the SCAD procedure works reasonably well, and the SCAD with a
correct EC working correlation structure provides noticeably better
variable selection results than the SCAD with IND working structure.

\section{Real data analysis}\label{sec:real}
In this section, the proposed methods are applied to analyze a
pharmacokinetics study for investigating CCI-779 effects on renal
cancer patients \cite{boni.leister.etc:2005}. CCI-779 is an anticancer
agent with demonstrated inhibitory effects on tumor growth. In this
study, patients with advanced renal cell carcinoma received CCI-779
treatment weekly until demonstrated evidence of disease progression.
One goal of the study is to identify transcripts in peripheral blood
mononuclear cells (PBMCs) which are useful for predicting the temporal
pharmacogenomic profile of CCI-799, after initiation of CCI-779
therapy. The data consists of expression levels of 12,626 genes from
$33$ patients on three scheduled visits: baseline, week 8 and week 16.
However, not all patients have measurements at all three visits. We
have unbalanced data with a total of only 54 observations. To account
for the cumulative-dose drug exposure, CCI-779 cumulative AUC was used
to quantify the pharmacogenomic measure of CCI-799 for each patient at
each visit. The AUC is of popular use in estimating bioavailability of
drugs in pharmacology. Since the response variable CCI-779 cumulative
AUC is continuous, we consider our model (\ref{mgaplm}) with an
identity link function.

With a total of 12,626 genes as covariates and only 54 observations, we
first apply the nonparametric independence screening method (NIS)
described in \cite{fan:feng:song:11} to reduce the dimensionality to a
moderate size. We ranked the genes according to their empirical
marginal function norms, and kept only the first $205$ genes with
marginal function norms larger than the 99th\% quantile of the
empirical norms of randomly permuted data. After variable screening, we
then applied the penalized polynomial splines \cite
{huang:horowitz:wei:10,xue:qu:12} for high-dimensional additive model
selection. We used the linear spline with a LASSO penalty function and
selected the tuning parameters with a five-fold cross-validation
procedure. This procedure further reduced dimensionality and selected
only 14 genes. Out of the selected 14 genes, we then applied our
proposed methods for more refined variable selection and estimation.

\begin{figure}

\includegraphics{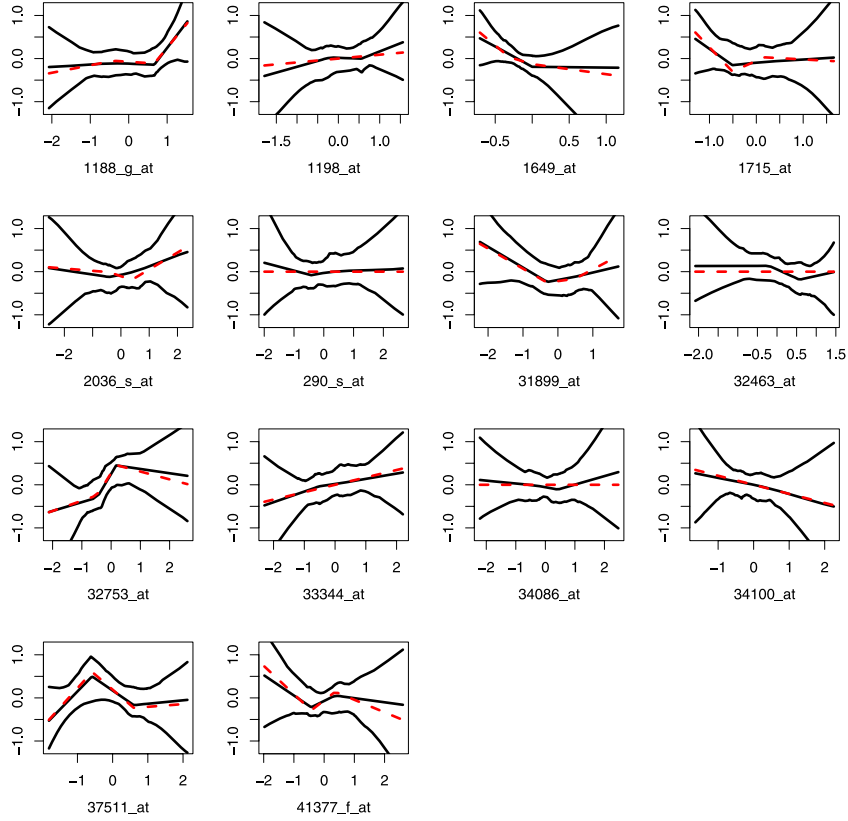}

\caption{Real data: Plots of the estimated components
of GAM (line) and GAPLM-SCAD (dashed) with $95$\% point-wise bootstrap
confidence intervals from the GAM.}
\label{realdata_plot}
\end{figure}

We first considered a generalized additive model (GAM), which is a
special case of a GAPLM model with $Z_{it}$ in (\ref{mgaplm})
consisting of an intercept term only. We applied the linear spline QIF
method to estimate the function components. The plots of the estimated
functions in Figure~\ref{realdata_plot} suggested that the function
forms of the five variables (1198\_at, 290\_s\_at, 32463\_at, 33344\_
at, 34100\_at) are almost linear. Therefore, we further considered a
GAPLM model with these five terms as linear terms, and the rest as
additive terms. For both models, we applied our proposed penalized QIF
method for more refined variable selection. For the GAPLM, we also
considered the variable selection method of \cite
{wang:liu:liang:carroll:11}. However, it can only select linear terms
and keeps all additive terms. We refer to this method as GAPLM-Linear.
Finally, as a benchmark, we also considered two linear models; one
contains only the $14$ genes selected in the high-dimensional additive
model and is referred as GLM, the other one begins with $205$ genes,
and variable selection in this high-dimensional linear model is then
conducted using LASSO, which is referred as GLM-LASSO.

For the GAM, we kept all $14$ variables, while both GAPLM and
GAPLM-Linear selected $11$ variables. In Table~\ref{Real_table}, we
report their mean squared estimation errors (MSEE) and EBIC values.
With the response being continuous, let $\hat{Y}_{it}=\sum_{l=1}^{d_x}\hat{\alpha}_{l} (x_{it}^{(l)} )+\mathbf
{z}_{it}^{t}\hat{\bolds{\beta}}$ be the estimator of $Y_{it}$ from
any method. Then define $\mbox{MSEE}=\frac{1}{N_{t}}\sum_{i=1}^{n}\sum_{t=1}^{T_i} (Y_{it}-\hat{Y}_{it} )^{2}$, with $N_{t}$ being
the total number of observations and $T_{i}$ being the size of cluster
$i$. Equation (\ref{EBIC}) with a Gaussian likelihood was used to
compute the EBIC, since the response variable is continuous and a
working independent structure is used here. It is not surprising that
the GAM gave the smallest MSEE since it has the most complicated model;
while the GAPLM-SCAD gives the most parsimonious model with the
smallest EBIC value. This suggests that with a simpler model, one may
be able to make more efficient estimation and inference. For the two
linear models, their much larger MSEEs suggest that the data contains
nonlinear dynamics which cannot be fully incorporated by linear models.

Furthermore, as suggested by one referee, we also compared the above
methods by their prediction performances. We randomly selected $28$
patients for estimation and left the remaining $5$ patients for
prediction. We calculated the mean squared prediction errors (MSPE) for
each method for $100$ replications. Table~\ref{Real_table} reports the
averaged MSPEs from $100$ replications. It shows that the GAPLM-SCAD
gives the smallest prediction error, and all non or semiparametric
methods give smaller prediction errors than the linear models. It again
suggests that the data contains a nonlinear structure. Those findings
are consistent with the ones observed from EBICs. In the above, we have
used an independent correlation structure in all procedures. Using
other types of correlation structure (e.g., exchangeable, AR-1) in the
estimation of GAM, GAPLM and GLM, which are not reported here, always
gives larger MSEEs due to the extra noise in modeling within-cluster
correlation when the sample size is rather small.

\begin{table}
\caption{Real data. The mean squared estimation error (MSEE), EBIC
values and averaged mean squared prediction error (AMSPE) of different methods}
\label{Real_table}
\begin{tabular*}{\textwidth}{@{\extracolsep{\fill}}lcccc@{}}
\hline
\textbf{Method} & \textbf{MSEE} & \textbf{EBIC} & \textbf{Model size} &\textbf{AMSPE} \\
\hline
GAM & $0.0127$ & \phantom{$-$}$0.1902$ & $14$& $0.4618$\\
GAPLM & $0.0162$ & $-0.6275$ & $14$& $0.3496$ \\
GAM-SCAD & $0.0132$ & \phantom{$-$}$0.2285$ & $14$& $0.2949$ \\
GAPLM-SCAD & $0.0205$ & $-0.8969$ & $11$&$0.2398$ \\
GAPLM-Linear &$0.0191$ &$-0.7774$ &$11$&$0.4069$\\
GLM & $0.2801$ & \phantom{$-$}$0.2772$ & $14$ &$0.6760$\\
GLM-LASSO & $0.0989$ & \phantom{$-$}$0.9716$ & $31$ &$0.8530$\\
\hline
\end{tabular*}
\end{table}

\section{Discussion}\label{sec:dics}

In this paper, we provide new statistical theory for model selection
and estimation with diverging numbers of linear covariates and
nonlinear components for generalized partial linear additive models.
Our work differs from existing works in three major aspects. First, we
consider model selection for both the parametric and nonparametric
parts simultaneously, while most of the literature focuses on selection
for either the parametric or the nonparametric part. Second, we allow
the numbers of linear covariates and nonlinear components to increase
with the sample size. Theoretical development for model selection and
estimation for diverging number of covariates in nonparametric
components is completely different from finite dimension settings.
Third, we allow dependence between the covariates in the nonparametric
and parametric part, and also dependence between the longitudinal
responses. All of these impose significant challenges in developing
asymptotic theory and oracle properties.

Note that the growing dimensions of the nonparametric part are smaller
than the parametric part, since the nonparametric components involve
many more parameters than the parametric part. The order of the
parametric dimension is comparable to that in the existing literature
for parametric model selection with diverging number of covariates
\cite{fan:peng:04,cho:qu:13,lian:liang:wang:13}. To establish the
asymptotic properties of the QIF estimators, a crucial step is to
obtain the upper and lower bounds of the eigenvalues of the matrix
$\mathbf{C}_{n}$ in the QIF equation. These bounds are assumed for
the parametric models \cite{fan:lv:08} or can be derived for
independent observations \cite{xue:09} using Lemma~1 of \cite
{stone:85}. However, neither of these are valid in our setting.
Instead, we develop an alternative strategy through proving Lemma~S.4,
which is essential in establishing bounds for the
eigenvalues of a large random matrix. The result in Lemma~S.4 can also be used for verifying the second-order KKT
condition on demand of bounds of random matrix with diverging dimension.

It is worth noting that the GEE estimator under the generalized partial
linear additive model framework is semiparametric efficient under the
correct correlation structure \mbox{\cite{cheng:zhou:huang:13}}. Since
the GEE and QIF are asymptotically equivalent when the correlation
structure is correctly specified, the proposed QIF estimator for the
generalized partial linear additive model is also semiparametric
efficient under the correct correlation structure.

\begin{appendix}
\section*{Appendix: Assumptions and proofs}\label{app}

\subsection{Notation and definitions}\label{subsec:app1}

For any functions $s_{1},s_{2}\in\mathcal{L}_{2} ( [0,1]
)$, define
%
\begin{equation}
\langle s_{1},s_{2} \rangle=\int_{0}^{1}s_{1}
( {x} ) s_{2} ( {x} ) \,d{x}\quad\mbox{and}\quad\Vert s \Vert_{2}^{2}=
\int_{0}^{1}s^{2} ( {x} ) \,d{x}.
\label{equ_inner}
\end{equation}
Let $\mathcal{H}_{0}$ be the space of constant functions
on $[0,1]$, and let $\mathcal{H}_{0}^{\bot}=\{s \dvtx \langle s,
1
\rangle=0,s \in\mathcal{L}_{2}\}$ and $1$ is the
constant function on $[0,1]$. Define the additive model space $%
\mathcal{M}$ and the space of additive polynomial spline functions $%
\mathcal{M}_{n}$ as
\[
\mathcal{M} = \Biggl\{ s ( \mathbf{x} ) =\sum_{l=1}^{d_{x}}s_{l}
( x_{l} );s_{l}\in \mathcal{H}_{0}^{\bot}
\Biggr\}, \qquad  \mathcal{M}_{n} = \Biggl\{ s ( \mathbf{x} ) =\sum
_{l=1}^{d_{x}}s_{l} (
x_{l} );s_{l}\in\varphi _{l}^{0,n}
\Biggr\},
\]
where\vspace*{1pt} $\varphi_{l}^{0,n}   =   \{ s_{l} ( \cdot
)
\dvtx s_{l}\in
\varphi_{l},  \langle s_{l},1 \rangle = 0 \}$ is the
centered polynomial spline space.
Let $\mathbf{s}  ( \mathbf{X} ) =  ( s
(
\mathbf{X}_{1} ),\ldots,s  (
\mathbf{X}_{T} )  )^{\T}$, for any $s\in\mathcal
{M}$ and
$\mathbf{X} = ( \mathbf{X}_{1}^{%
\T},\ldots,\mathbf{X}_{T}^{\T} ) ^{\T}$. We define the
theoretical and empirical norms
of $\mathbf{s}$: $\Vert\mathbf{s}\Vert^{2}=E \{
\mathbf
{s}^{\T} ( \mathbf{X%
} ) \mathbf{s}  ( \mathbf{X} )  \} $ and
$\Vert\mathbf{s} \Vert
_{n}^{2}=\frac{1}{n}\sum_{i=1}^{n}\mathbf{s}^{\T} (
\mathbf
{X}_{i} )
\mathbf{s} ( \mathbf{X}_{i} ) $.

For $\bolds{\theta}=(\bolds{\beta}^{\T},\bolds{\gamma
}^{\T})^{\T}$,
denote a $d_{n}$ vector and a $d_{n}\times d_{n}$ matrix
%
\begin{equation}\qquad
\mathbf{S}_{n} ( \bolds{\theta} )=\dot{\mathbf
{G}}_{n}^{\T%
} ( \bolds{\theta} ) \mathbf{C}_{n}^{-1}
( \bolds{\theta}%
) \mathbf{G}_{n} ( \bolds{
\theta} ),  \qquad  \mathbf{H}_{n} ( \bolds{\theta} )=\mathbf {
\dot {G}}_{n}^{\T%
} ( \bolds{\theta} )
\mathbf{C}_{n}^{-1} ( \bolds{\theta}%
)
\dot{\mathbf{G}}_{n} ( \bolds{\theta} ), \label{DEF:SnHn}
\end{equation}
where the $(Kd_{n})\times d_{n}$ matrix
%
\begin{equation}\qquad
\dot{\mathbf{G}}_{n} ( \bolds{\theta} ) \equiv
\frac
{\partial}{%
\partial\bolds{\theta}}\mathbf{G}_{n} ( \bolds{\theta} ) = \biggl(
\frac{\partial}{\partial\bolds{\beta}}\mathbf {G}_{n} ( \bolds{\theta} ),
\frac{\partial}{\partial\bolds{\gamma}}\mathbf{G%
}_{n} ( \bolds{\theta}
) \biggr) \equiv \bigl( \dot{\mathbf{G}}_{%
\bolds{\beta}} ( \bolds{
\theta} ),\dot{\mathbf {G}}_{\bolds{%
\gamma}} ( \bolds{\theta} ) \bigr).
\label{EQ:Gdot}
\end{equation}
By \cite{qu:lindsay:li:00} and Lemma~S.4, the estimating
equation for $\bolds{\theta}$ is
%
\begin{equation}
n^{-1}\dot{Q}_{n} ( \bolds{\theta} ) \equiv
n^{-1}\frac
{\partial
}{\partial\bolds{\theta}}Q_{n} ( \bolds{\theta } ) =2
\mathbf{S}%
_{n} ( \bolds{\theta} ) +
O_{P}{ \bigl(%
\rho_n^{-1}n^{-1}
d_{n} \bigr)} =\mathbf{0}, \label{EQ:Qdiv1}
\end{equation}
and the second derivative of $Q_{n} ( \bolds{\theta}
) $
in $%
\bolds{\theta}$
%
\begin{equation}
n^{-1}\ddot{Q}_{n} ( \bolds{\theta} ) \equiv
n^{-1}\frac
{\partial
^{2}}{\partial\bolds{\theta}\,\partial\bolds{\theta}^{\T
}}Q_{n} ( \bolds{\theta} ) =2
\mathbf{H}_{n} ( \bolds{\theta} ) +o_{P}(1).
\label{EQ:Qdiv2}
\end{equation}
To facilitate technical arguments in the following proofs, we write
%
\begin{eqnarray}
\mathbf{S}_{n}(\bolds{\theta}) 
&=& \pmatrix{ \dot{\mathbf{G}}_{\bolds{\beta}}^{\T}
( \bolds{\theta} ) \mathbf{C}_{n}^{-1} (
\bolds{\theta} ) \mathbf {G}_{n} ( \bolds{\theta} )
\vspace*{2pt}\cr
\dot{\mathbf{G}}_{\bolds{\gamma}}^{\T} ( \bolds{\theta} )
\mathbf{C}_{n}^{-1} ( \bolds{\theta} ) \mathbf
{G}_{n} ( \bolds{\theta} )}, \label{EQ:Sn}
\\
%
\mathbf{H}_{n}(\bolds{\theta}) 
&=& \pmatrix{ \mathbf{
\dot{G}}_{\bolds{\beta}}^{\T} ( \bolds{\theta} )
\mathbf{C}_{n}^{-1} ( \bolds{\theta} ) \mathbf {
\dot{G}}_{\bolds{%
\beta}} ( \bolds{\theta} ) & \dot{\mathbf
{G}}_{\bolds{\beta}}^{%
\T} ( \bolds{\theta} ) \mathbf{C}_{n}^{-1}
( \bolds{\theta}%
) \dot{\mathbf{G}}_{\bolds{\gamma}} (
\bolds{\theta} )
\vspace*{2pt}\cr
\dot{\mathbf{G}}_{\bolds{\gamma}}^{\T} ( \bolds{\theta} )
\mathbf{C}_{n}^{-1} ( \bolds{\theta} ) \mathbf {
\dot{G}}_{\bolds{%
\beta}} ( \bolds{\theta} ) & \dot{\mathbf
{G}}_{\bolds{\gamma}%
}^{\T} ( \bolds{\theta} ) \mathbf
{C}_{n}^{-1} ( \bolds{%
\theta} )
\dot{\mathbf{G}}_{\bolds{\gamma}} ( \bolds{\theta}%
)%
}. \label{EQ:Hn}
\end{eqnarray}

\subsection{Assumptions}
\label{subsec:assu}

We denote $ ( \mathbf{Y}_{i},\mathbf{X}_{i},\mathbf
{Z}_{i} )
,i=1,\ldots,n$ which are i.i.d. samples from population $ (
\mathbf{Y},%
\mathbf{X},\mathbf{Z} ) $ with $\mathbf{Y}= (
Y_{1},\ldots
,Y_{T} ) ^{\T}$, $\mathbf{X}= ( \mathbf
{X}_{1},\ldots
,\break \mathbf{X}%
_{T} ) ^{\T}$, and $\mathbf{Z}= ( \mathbf
{Z}_{1},\ldots,\mathbf{Z}%
_{T} ) ^{\T}$ for correlated data with cluster size $T$. Denote
$C^{(r)} ( [0,1] )=\{
\mbox{$f$: $f$ has continuous derivatives up to
order $r$ on } [0,1]\}$
as the space of the $r$th order smooth functions on $[0,1]$.
For any vector $\mathbf{a}$, let $\Vert
\mathbf{a}\Vert$
be the usual Euclidean norm. For any matrix $\mathbf{A}$, let
$\Vert
\mathbf{A}\Vert$ be the modulus of the largest singular value of
$\mathbf{A}$. We provide the regularity conditions to obtain
Theorems \ref{THM:nonparametric}--\ref{THM:oracle}.

\begin{longlist}[(C4)]
\item[(C1)] For some $r\geq2$, $\alpha_{0,l}\in C^{(r)} (
[0,1] )$
$l=1,\ldots,d.$

\item[(C2)] \textit{The covariance matrix} $\bolds{\Sigma
}=E\mathbf{ee}^{\T%
}$ \textit{is positive definite, and} $E\Vert\mathbf
{e}\Vert
^{2+\delta}<+\infty$ \textit{for some} $\delta>0$.

\item[(C3)] \textit{For each $\mathbf{X}_{t}$}, $t=1,\ldots,T$,
its density
function $f_{t} ( \mathbf{x} )$ \textit{is
absolutely continuous
and bounded away from zero and $\infty$ on a compact support $\chi
=[0,1]^{d_x}$.}

\item[(C4)] \textit{The number of nonzero components in the nonparametric
part} $s_{x}$ \textit{is fixed; there exists} $c_{\alpha}> 0$
\textit{such that} $\min_{1\leq l \leq s_x} \|\alpha_{0,l}\|>
c_{\alpha
}$.
\textit{The nonzero coefficients in the linear part satisfy that} $%
\min_{1\leq k \leq s_z} \|\beta_{0k}\|/\lambda_{2n}\rightarrow
\infty$.

\item[(C5)] \textit{The eigenvalues
of} $E( \bolds{\Gamma}_{0}^{(k)}) $ \textit{are bounded away
from} $0$ \textit{and $\infty$, uniformly in} $k=1,\ldots,K$,
\textit{for sufficiently large} $n$.

\item[(C6)] \textit{The second derivative of}
$g^{-1} ( \cdot ) $ \textit{exists and is bounded; function
$V (\cdot )$
has a bounded second derivative, and is bounded away from $0$ and
$\infty$}.

\item[(C7)] 
\textit{The modular of the singular
value of} $\mathbf{M}= ( \mathbf{M}_{1}^{\T},\ldots,%
\mathbf{M}_{K}^{\T} ) ^{\T}$ \textit{is bounded away from} $0$
\textit{and $\infty$}.

\item[(C8)] \textit{The eigenvalues of} $E(\mathbf
{X}_{t}\mathbf
{X} _{t}^{\T%
}\vert\mathbf{Z}_{t})$ \textit{are bounded away from}
$0$ \textit{and $\infty$, uniformly in} $1\leq t\leq T$.

\item[(C9)]
{\textit{There is a large enough open subset} $\widetilde{\Theta
}_n\in R^{d_n}$
\textit{which contains} $\widetilde{\bolds{\theta
}}_{0}=(\bolds{\beta}_{0}^{\T},
\widetilde{\bolds{\gamma}}^{\T})^{\T}$, \textit{for}
$\widetilde
{\bolds{\gamma}}$
\textit{in Section~\ref{subsec:nonparametric}, such that} $\sup_{\bolds{\theta}\in\widetilde{\Theta}_n}
\llvert n^{-1}\frac{\partial^3 Q_{n} ( \bolds{\theta
} )
}{\partial\theta_{j}\,\partial\theta_{k} \,\partial\theta_{l}}\rrvert =O_{P}(\rho_n^{-1})$.}

\item[(P1)] $\liminf_{n \rightarrow\infty}\liminf_{\theta
\rightarrow0+}
p_{\lambda_{j,n}}(\theta)/\lambda_{j,n}>0$, $j=1,2$.

\item[(P2)] $a_n=o(1/\sqrt{nd_n})$, $b_n=o(d_n^{-1/2})$,
\textit{where} $a_n$
\textit{and} $b_n$ \textit{are defined in~(\ref{DEF:anbn}).}
\end{longlist}

Conditions (C1)--(C3) are quite standard in the spline smoothing
literature. Assumptions similar to (C1)--(C3) can be found in \cite
{huang:98,huang:03,xue:09} and \cite{xue:qu:zhou:10}. The
smoothness condition in (C1) controls the rate of convergence of the
spline estimators $ \widehat{\alpha}_{l}$, \textcolor
{black}{$l=1,\ldots,d_{x}$}, and $\widehat{ \alpha}$. Conditions (C5) and
(C6) are similar to assumptions (A3) and (A4) in \cite{he:fun:zhu:05},
which can be verified for other distributions as well. The boundedness
condition in condition (C7) is essentially a requirement that the
matrix $\mathbf{C}_{n}$ in (\ref{EQ:alphahat}) is asymptotically
positive definite. This assumption is clearly satisfied if the basis
matrices are exchangeable or AR-1 correlation structures as discussed
previously. The condition on eigenvalues in (C8) is to ensure that we
do not have a multicolinear problem. Condition (C9) controls the
magnitude of the third-order derivative of the quadratic inference
function. Similar conditions have been assumed in \cite{cho:qu:13} and
\cite{fan:peng:04}. Here, we require a slightly stronger condition.
Instead of assuming boundedness, we require it be of the order
$O_{P}(\rho_{n}^{-1})$, where $\rho_n=((1-\delta)/2)^{(d_x-1)/2}$ to
facilitate the technical derivation for the nonparametric components in
a GAPLM model, while both \cite{cho:qu:13} and \cite{fan:peng:04}
consider pure parametric models.

\subsection{Proof of Theorem \texorpdfstring{\protect\ref{THM:nonparametric}}{1}}
\label{subsec:nonparametric}

According to Lemma A.7 of \cite{xue:qu:zhou:10}, for any function $%
\alpha\in\mathcal{M}$ with $\alpha_{l}\in C^{(r)} (
[0,1]
)$, $%
l=1,\ldots,d_x$, there exists an additive spline function $\widetilde
{\alpha}=%
\widetilde{\bolds{\gamma}}^{\T}\mathbf{B}\in\mathcal{M}_{n}$
and a
constant $C$ such that
%
\begin{equation}
\Vert\widetilde{\alpha}-\alpha\Vert_{\infty}\leq Cd_{x}J_n^{-r}.
\label{DEF:alphatilde}
\end{equation}
From the results of Lemma~S.10 in the online
supplementary material \cite{wang:xue:qu:liang:13} and Lemma A.6 in
the online supplement of \cite{xue:qu:zhou:10}, we have
%
\begin{eqnarray}\label{equ_Bnorm}
\bigl\Vert\mathbf{B}^{\T}\bigl(\widehat{\bolds{
\gamma}}^{\qif
}-\widetilde{\bolds{%
\gamma}}\bigr)
\bigr\Vert_{n}^{2} &=&\frac{1}{n}\sum
_{l=1}^{d_x}\sum_{i=1}^{n}%
\sum_{t=1}^{T} \bigl[
\mathbf{B}_{l}^{\T} \bigl(x_{it}^{(l)}
\bigr) \bigl(%
\widehat{\bolds{\gamma}}_{l}^{\qif}-
\widetilde{\bolds{\gamma }}_{l}\bigr)%
\bigr]
^{2}
\nonumber
\\[-8pt]
\\[-8pt]
\nonumber
& =&O_{P}\bigl(n^{-1} d_n\bigr).
\end{eqnarray}
The triangular inequality implies that, for each $l=1,\ldots,d_{x}$,
\begin{eqnarray*}
\frac{1}{n}\sum_{l=1}^{d_x}\sum
_{i=1}^{n}\sum
_{t=1}^{T} \bigl[ \widehat{%
\alpha}_{l} \bigl(x_{it}^{(l)}\bigr) -
\alpha_{0,l}\bigl( x_{it}^{(l)}\bigr) \bigr]
^{2} & \leq & \frac{2}{n}\sum_{l=1}^{d_x}
\sum_{i=1}^{n} \sum
_{t=1}^{T} \bigl[ \mathbf{B}%
_{l}^{\T} \bigl(x_{it}^{(l)}\bigr) \bigl(
\widehat{\bolds{\gamma}}_{l}^{\qif
}-\widetilde{%
\bolds{\gamma}}_{l}\bigr) \bigr] ^{2}
\\
&&{}+Cd_{x}J_n^{-2r}.
\end{eqnarray*}
This completes the proof.

\subsection{Proof of Theorem \texorpdfstring{\protect\ref{THM:parametric}}{2}}

To study the asymptotic properties of $\widehat{\bolds{\beta
}}^{\qif}$, we consider the case that $\alpha_0$ in (\ref{mgaplm}) can
be estimated at reasonable
accuracy, for example, we can approximate $\alpha_0$ by the spline\vadjust{\goodbreak}
smoother $%
\widetilde{\alpha}$ in (\ref{DEF:alphatilde}). We begin our proof by
replacing $\alpha_0$ with $\widetilde{\alpha}$ and
defining an \textit{intermediate} QIF estimator for $\bolds{\beta}_0$.

For any fixed $\bolds{\beta}$ and $i=1,\ldots,n$, we denote
$\widetilde{%
\bolds{\eta}}_{i}(\bolds{\beta})=\widetilde{\alpha
}(\mathbf{X}_{i})+%
\mathbf{Z}_{it}^{\T}\bolds{\beta}$ and $\widetilde
{\bolds{\mu}}_{i}(%
\bolds{\beta})=g^{-1}\{\widetilde{\bolds{\eta
}}_{i}(\bolds{\beta})\}$.
Let $\dot{\widetilde{\mu}}_{it}(%
\bolds{\beta})$ be the first-order derivative\vspace*{1pt} of $g^{-1}(\eta)$
evaluated at $%
\eta=\widetilde{\bolds{\eta}}_{i}(\bolds{\beta})$. Define
$\widetilde{\bolds{\Delta}}_{i}(\bolds{\beta})=\diag%
\{ \dot{\widetilde{\mu}}_{i1}(\bolds{\beta}),\ldots
,\dot
{\widetilde{%
\mu}}_{iT}(\bolds{\beta}) \} $ and $\widetilde
{\mathbf
{A}}%
_{i}=\break \diag \{ V(\widetilde{\mu}_{i1}),\ldots, V(\widetilde{\mu
}_{iT}) \} $. Let
%
\begin{equation}
\widetilde{\mathbf{g}}_{i}(\bolds{\beta})=\mathbf
{g}_{i} ( \bolds{%
\beta,}\widetilde{\bolds{
\gamma}} ) = \pmatrix{ \mathbf{D}_{i}^{\T}
\widetilde{\bolds{\Delta }}_{i}\widetilde {\mathbf{A}}%
_{i}^{-1/2}\mathbf{M}_{1}\widetilde{
\mathbf {A}}_{i}^{-1/2} \bigl( \mathbf{Y}%
_{i}-\widetilde{\bolds{\mu}}_{i}(\bolds{\beta})
\bigr)
\vspace*{2pt}\cr
\vdots
\vspace*{2pt}\cr
\mathbf{D}_{i}^{\T}\widetilde{\bolds{\Delta
}}_{i}\widetilde {\mathbf{A}}%
_{i}^{-1/2}
\mathbf{M}_{K}\widetilde{\mathbf {A}}_{i}^{-1/2}
\bigl( \mathbf{Y}%
_{i}-\widetilde{\bolds{
\mu}}_{i}(\bolds{\beta}) \bigr)}. \label{DEF:gitilde}
\end{equation}
Define $\widetilde{\mathbf{G}}_{n}(\bolds{\beta})=\frac
{1}{n}\sum_{i=1}^{n}%
\widetilde{\mathbf{g}}_{i} ( \bolds{\beta} ) $.
In a similar
way, we define $\widetilde{\mathbf{C}}_{n}(\bolds{\beta})$,
and $%
\widetilde{Q}_{n}(\bolds{\beta})$. Let $\widetilde{\bolds{\beta}}_{\qif}=
\argmin_{\bolds{\beta}}n^{-1}\widetilde{Q}_{n}(\bolds{\beta
})=\argmin_{\bolds{\beta}} \{\widetilde{\mathbf
{G}}_{n}^{\T}(%
\bolds{\beta})\widetilde{\mathbf{C}}_{n}^{-1}(\bolds{\beta})\widetilde{%
\mathbf{G}}_{n}(\bolds{\beta}) \}$.
The\break  asymptotic properties of $\widetilde{\bolds{\beta}}_{\qif}$
are given
in the supplementary material \cite{wang:xue:qu:liang:13}. Let $%
\widehat{\bolds{\theta}}^{\qif}=(\widehat{\bolds{\beta
}}_{\qif}^{\T},%
\widehat{\bolds{\gamma}}_{\qif}^{\T})^{\T}$, $\widetilde
{\bolds{\theta}}%
_{0}=(\bolds{\beta}_{0}^{\T},\widetilde{\bolds{\gamma
}}^{\T
})^{\T}$ and $%
\widetilde{\bolds{\theta}}^{\qif}=(\widetilde{\bolds{\beta
}}_{\qif}^{\T},%
\widetilde{\bolds{\gamma}}^{\T})^{\T}$.

\begin{pf*}{Proof of Theorem~\ref{THM:parametric}} By Taylor expansion,
\begin{eqnarray*}
&&\dot{Q}_{n}\bigl(\widehat{\bolds{\theta}}^{\qif}\bigr)-
\dot {Q}_{n}(\widetilde{%
\bolds{\theta}}_{0})\\
&&\qquad=
\ddot{Q}_{n}( \widetilde{\bolds{\theta}}_{0}) \bigl(
\widehat{\bolds{\theta}}^{\qif}-\widetilde{\bolds{\theta
}}_{0}\bigr)+\frac{1}{%
2}\bigl(\widehat{\bolds{
\theta}}^{\qif}-\widetilde{\bolds{\theta }}_0
\bigr)^{\T%
}\frac{\partial\dot{Q}_{n} ( \bolds{\theta} )
}{\partial
\bolds{\theta}\,\partial\bolds{\theta}^{\T}}\bigg\rrvert _{\bolds{\theta}=%
\bolds{\theta}^{\ast}} \bigl(
\widehat{\bolds{\theta }}^{\qif} -\widetilde{%
\bolds{
\theta}}_0\bigr),
\end{eqnarray*}
where $\bolds{\theta}^{\ast}=t\widehat{\bolds{\theta
}}^{\qif}+(1-t)%
\widetilde{\bolds{\theta}}_{0}$, for some $t\in{}[0,1]$.
Since $%
\dot{Q}_{n}(\widehat{\bolds{\theta}}^{\qif})=0$,
\[
-\dot{Q}_{n}(\widetilde{\bolds{\theta}}_{0})=
\ddot{Q}_{n}(%
\widetilde{\bolds{\theta}}_{0})
\bigl(\widehat{\bolds{\theta }}^{\qif}-%
\widetilde{
\bolds{\theta}}_{0}\bigr)+\frac{1}{2}\bigl(\widehat {
\bolds{\theta}%
}^{\qif}-\widetilde{\bolds{
\theta}}_0\bigr)^{\T}\frac
{\partial
\dot{Q}%
_{n} ( \bolds{\theta} ) }{\partial\bolds{\theta
}\,\partial
\bolds{\theta}^{\T}}\bigg\rrvert
_{\bolds{\theta}=\bolds{\theta}^{\ast}}\bigl(%
\widehat{\bolds{\theta}}^{\qif}-
\widetilde{\bolds{\theta}}_0\bigr).
\]
According to the Cauchy--Schwarz inequality, one has
\[
\biggl\llVert \frac{1}{n}\bigl(\widehat{\bolds{\theta}}^{\qif
}-
\widetilde {\bolds{\theta}}_0\bigr)^{\T}
\frac{\partial\dot{Q}_{n} ( \bolds{\theta} )
}{\partial
\bolds{\theta}\,\partial\bolds{\theta}^{\T}} \bigl(\widehat{\bolds{\theta}}^{\qif}-\widetilde{
\bolds{\theta }}_0\bigr)\biggr\rrVert ^{2} \leq\bigl\|
\widehat{\bolds{\theta}}^{\qif}-\widetilde {\bolds{
\theta}}_0\bigr\|^{4} \frac{1}{n} \sum
_{j,k,l=1}^{d_n} \biggl\{\frac{\partial^3
Q_{n} (
\bolds{\theta} ) } {
\partial\theta_{j}\,\partial\theta_{k}\, \partial\theta_{l}} \biggr
\}^{2}.
\]
Lemma~S.10 and condition (C9) implies that
\begin{eqnarray*}
\biggl\llVert \frac{1}{n}\bigl(\widehat{\bolds{\theta}}^{\qif
}-
\widetilde {\bolds{\theta}}_0\bigr)^{\T}
\frac{\partial\dot{Q}_{n} ( \bolds{\theta} )
}{\partial
\bolds{\theta}\,\partial\bolds{\theta}^{\T}} \bigl(\widehat{\bolds{\theta}}^{\qif}-\widetilde{
\bolds{\theta }}_0\bigr)\biggr\rrVert ^{2} & \leq &
\bigl\|\widehat{\bolds{\theta}}^{\qif
}-\widetilde {\bolds{
\theta}}_0\bigr\|^{4}\times O_{P}\bigl(
\rho_n^{-1}d_n^3\bigr)
\\
&=&O_{P}\bigl(n^{-2} d_n^2
\bigr)\times O_{P}\bigl(\rho _n^{-1}d_n^3
\bigr)=o_{P}\bigl(n^{-1}\bigr).
\end{eqnarray*}
Next by (\ref{EQ:Qdiv1}) and (\ref{EQ:Qdiv2}), we have
\begin{eqnarray*}
- \bigl\{ 2\mathbf{S}_{n}(\widetilde{\bolds{
\theta}}_{0}) + O_{P} \bigl(%
{
\rho_n^{-1}}n^{-1}d_n \bigr) \bigr
\}& = & \bigl\{ 2\mathbf{H}%
_{n}(\widetilde{\bolds{
\theta}}_{0}) +o_{P}(1) \bigr\} \bigl(\widehat {
\bolds{%
\theta}}^{\qif} - \widetilde{\bolds{
\theta}}_{0}\bigr)
\\
&&{} +o_P\bigl(n^{-1/2}\bigr),
\end{eqnarray*}
where $\mathbf{S}_{n} ( \bolds{\theta} ) $ and
$\mathbf{H}%
_{n} ( \bolds{\theta} ) $ are defined in (\ref
{DEF:SnHn}). Thus,
\begin{eqnarray*}
\pmatrix{ \widehat{\bolds{\beta}}^{\qif}-
\bolds{\beta}_{0}
\vspace*{2pt}\cr
\widehat{\bolds{\gamma}}^{\qif}-\widetilde{\bolds{
\gamma}}} & = & - \left[2 \pmatrix{ \mathbf{H}_{\bolds{\beta\beta}} & \mathbf
{H}_{\bolds{\beta\gamma}}
\vspace*{2pt}\cr
\mathbf{H}_{\bolds{\gamma\beta}} & \mathbf {H}_{\bolds{\gamma\gamma}}%
} +o_{P}(1) \right]^{-1}
\\
&&{}\times \bigl[2\mathbf{S}_{n}(\widetilde{\bolds{\theta
}}_{0})+O_{P} \bigl({ \rho_{n}^{-1}}n^{-1}
d_n \bigr) \bigr]+o_P\bigl(n^{-1/2}\bigr),
\end{eqnarray*}
which leads to
\begin{eqnarray*}
\widehat{\bolds{\beta}}^{\qif}-\bolds{\beta}_{0}& =
& \bigl\{ \mathbf{H}_{\bolds{\beta\beta}}-\mathbf {H}_{\bolds{\beta\gamma}}
\mathbf{H}_{%
\bolds{\gamma\gamma}}^{-1}\mathbf{H}_{\bolds{\gamma
\beta}} \bigr
\} ^{-1} \bigl( \mathbf{I,H}_{\bolds{\beta\gamma}}\mathbf
{H}_{\bolds{\gamma
\gamma}}^{-1} \bigr) \mathbf{S}_{n}(
\widetilde{\bolds{\theta }}%
_{0})
\\
&&{}+O_{P} \bigl({ \rho_{n}^{-1}}n^{-1}d_n
\bigr)+o_P\bigl(n^{-1/2}\bigr).
\end{eqnarray*}
According to (\ref{EQ:Sn}),
\begin{eqnarray*}
&& \bigl( \mathbf{I,H}_{\bolds{\beta\gamma}} ( \bolds{\theta} )
\mathbf{H}_{\bolds{\gamma\gamma}}^{-1} ( \bolds{\theta} ) \bigr)
\mathbf{S}_{n} ( \bolds{\theta} )
\\
&&\qquad = \bigl( \mathbf{I,H}_{
\bolds{\beta\gamma}} ( \bolds{\theta} )
\mathbf{H}_{\bolds{%
\gamma\gamma}}^{-1} ( \bolds{\theta} ) \bigr)
\bigl( \mathbf{%
\dot{G}}_{\bolds{\beta}} ( \bolds{\theta} ),
\dot{\mathbf{G}}_{%
\bolds{\gamma}} ( \bolds{\theta} ) \bigr) ^{\T
}
\mathbf{C}%
_{n}^{-1} ( \bolds{\theta} )
\mathbf{G}_{n} ( \bolds{\theta}%
)
\\
&&\qquad = \bigl\{ \dot{\mathbf{G}}_{\bolds{\beta
}}^{\T
} ( \bolds{
\theta}%
) -\mathbf{H}_{\bolds{\beta\gamma}} ( \bolds{\theta} )
\mathbf{H}_{\bolds{\gamma\gamma}}^{-1} ( \bolds{\theta} )
\dot{\mathbf{G}}_{\bolds{\gamma}}^{\T} ( \bolds{\theta} ) \bigr
\} \mathbf{C}_{n}^{-1} ( \bolds{\theta} )
\mathbf{G}%
_{n} ( \bolds{\theta} ).
\end{eqnarray*}
Hence, the asymptotic distribution of $\sqrt{n}\mathbf
{A}_n\bolds{\Sigma}%
_{n}^{-1/2}(\widehat{\bolds{\beta}}^{\qif}-\bolds{\beta
}_{0})$ is the
same as that of
\begin{eqnarray*}
&&\sqrt{n}\mathbf{A}_n\bolds{\Sigma}_{n}^{-1/2}
\bigl\{ \mathbf {H}_{\bolds{\beta
\beta}}(\widetilde{\bolds{
\theta}}_{0})-\mathbf {H}_{\bolds{\beta\gamma}%
}(\widetilde{\bolds{
\theta}}_{0})\mathbf{H}_{\bolds{\gamma\gamma}}^{-1}(%
\widetilde{\bolds{\theta}}_{0})\mathbf{H}_{\bolds{\gamma
\beta}}(%
\widetilde{\bolds{\theta}}_{0})\bigr\}^{-1}
\\
&&\qquad\times \bigl\{ \bigl( \dot{\mathbf {G}}_{\bolds{\beta}}^{\T}(%
\widetilde{\bolds{\theta}}_{0})-\mathbf{H}_{\bolds{\beta
\gamma}}
\mathbf{%
H}_{\bolds{\gamma\gamma}}^{-1}\dot{\mathbf
{G}}_{\bolds{\gamma}}^{\T}(%
\widetilde{\bolds{
\theta}}_{0}) \bigr) \mathbf {C}_{n}^{-1}(
\widetilde{%
\bolds{\theta}}_{0})\mathbf{G}_{n}(
\widetilde{\bolds{\theta}}%
_{0}) \bigr\}.
\end{eqnarray*}
The desired result follows from Lemmas S.11 and~S.12. 
\end{pf*}
\subsection{Proof of Theorem \texorpdfstring{\protect\ref{THM:SCAD-QIF}}{3}}

In the following, let $L_{n}(\bolds{\theta})=Q_{n}(\bolds{\theta})+\break  n\sum_{l=1}^{d_{x}}p_{%
\lambda_{1,n}} ( \Vert\bolds{\gamma}_{l}\Vert
_{\mathbf
{K}%
_{l}} ) +\sum_{l=1}^{d_{z}}p_{\lambda_{2,n}}(|\bolds{\beta}_{l}|)$
be the object function in (\ref{DEF:pqif}). Let $\Theta_{\mathcal
{A}}=\{%
\bolds{\theta}=(\bolds{\beta}^{\T},\bolds{\gamma
}^{\T
})^{\T}\dvtx \beta
_{s_{z}+1}=\cdots=\beta_{d_{z}}=0, \bolds{\gamma
}_{s_{x}+1}=\cdots=%
\bolds{\gamma}_{d_{x}}=\mathbf{0}\}$ and define
$\widehat{\bolds{\theta}}_{\mathcal{A}}=(\widehat{\bolds{\beta}}_{%
\mathcal{A}}^{\T},\widehat{\bolds{\gamma}}_{\mathcal{A}}^{\T
})^{\T
}=%
\argmin_{\bolds{\theta}\in\Theta_{\mathcal
{A}}}Q_{n}(\bolds{\theta})$,
which leads to the spline QIF estimator of the nonzero components, given
that the rest terms are zero. Note that
$\Vert\mathbf{B}^{\T} ( \widehat{{\bolds{\gamma
}}}_{\mathcal{A}}-%
\widetilde{\bolds{\gamma}} ) \Vert
_{n}=O_{P}(n^{-1/2}d_n^{1/2})$ and
$\Vert\widehat{\bolds{\beta}}_{\mathcal{A}}-\bolds{\beta
}_{0}\Vert
=O_{P}(n^{-1/2}d_n^{1/2})$ from the results of Theorems \ref{THM:nonparametric}
and \ref{THM:parametric}. It is
sufficient to show that for large $n$ and any $%
\epsilon>0$, there exists a sufficient large constant $C$ such that
%
\begin{equation}
P \Bigl\{\inf_{\Vert\bolds{\theta}-\widehat{\bolds{\theta
}}_{\mathcal{A}%
}\Vert= C { \rho_{n}^{-3}}d_n^{1/2}(n^{-1/2}+a_n)} L_{n} ( \bolds{\theta} )
>L_{n}(\widehat{\bolds{\theta}}_{%
\mathcal{A}}) \Bigr\} \geq1-
\epsilon. \label{equ_th2}
\end{equation}
Equation (\ref{equ_th2}) implies that $L_{n} ( \cdot ) $ has
a local
minimum in the set
$ 
\Theta^{*}(C)=\{ \bolds{\theta}\dvtx \Vert\bolds{\theta
}-\widehat{\bolds{%
\theta}}_{\mathcal{A}}\Vert\leq C { \rho_{n}^{-3}}
d_n^{1/2}(n^{-1/2}+a_n)\}$. 
Thus, one has $\Vert\widehat{\bolds{\theta}}^{\qif}-\widehat
{\bolds{%
\theta}}_{\mathcal{A}}\Vert=\break  O_{P}\{{ \rho_{n}^{-3}}%
d_n^{1/2}(n^{-1/2}+a_n)\}$. Further, the triangular inequality yields
that\break
$\Vert\widehat{\bolds{\theta}}^{\qif}-\bolds{\theta
}_{0}\Vert\leq
\Vert\widehat{\bolds{\theta}}^{\qif}-\widehat{{\bolds{\theta}}}_{%
\mathcal{A}}\Vert+\Vert\widehat{{\bolds{\theta}}}_{\mathcal
{A}}-\bolds{%
\theta}_{0}\Vert=O_{P}\{{ \rho_{n}^{-3}}
d_n^{1/2}(n^{-1/2}+a_n)\}$. The theorem follows from condition (C4).

In the following, we show that (\ref{equ_th2}) holds. Observing that
$p_{\lambda_{n}} ( 0 ) =0$ and $p_{\lambda_{n}} (
\cdot
) \geq0$, one has
\begin{eqnarray*}
L_{n}(\bolds{\theta}) - L_{n}(
\widehat {\bolds{\theta}}_{\mathcal{A}}) &\geq &Q_{n}(\bolds{
\theta}) - Q_{n}(\widehat{\bolds{\theta }}_{\mathcal{A}%
}) +
\sum_{l=1}^{s_{x}}n \bigl\{ p_{\lambda_{1,n}} \bigl(
\Vert \bolds{\gamma}%
_{l}\Vert_{\mathbf{K}_{l}} \bigr) -
p_{\lambda_{1,n}}\bigl ( \Vert\widehat{{%
\bolds{\gamma}}}_{\mathcal{A}}
\Vert_{\mathbf
{K}_{l}} \bigr) \bigr\}
\\
&&{} +\sum_{l=1}^{s_{z}}n\bigl\{
p_{\lambda_{2n}} \bigl( |\beta _{l}| \bigr) -p_{\lambda_{2,n}}\bigl(\vert\widehat{
\beta}_{\mathcal
{A},l}\vert \bigr)\bigr\}.
\end{eqnarray*}
Note that
\begin{eqnarray*}
Q_{n}(\bolds{\theta}) - Q_{n}(\widehat{\bolds{
\theta }}_{\mathcal{A}})=(%
\bolds{\theta}-\widehat{\bolds{
\theta}}_{\mathcal
{A}})^{\T
}\dot{Q}_{n}(%
\widehat{\bolds{\theta}}_{\mathcal{A}}) + \tfrac
{1}{2}(\bolds{
\theta}-%
\widehat{\bolds{\theta}}_{\mathcal{A}})^{\T}
\ddot {Q}_{n}(\widehat {\bolds{%
\theta}}_{\mathcal{A}})
(\bolds{\theta} - \widehat {\bolds{\theta}}_{%
\mathcal{A}})+R_{n}^{\ast},
\end{eqnarray*}
where
\[
\bigl|R_{n}^{\ast}\bigr|=\frac{C^3}{6}\Biggl\llvert \sum
_{i,j,k=1}^{d_n}\frac
{\partial
Q_{n}(\bolds{\theta}^{*})}{\partial\theta_{i}\, \partial\theta_{j}
\,\partial
\theta_{k}}(\theta_{i}-
\widehat{\theta}_{\mathcal{A},i}) (\theta _{j}-%
\widehat{
\theta}_{\mathcal{A},j}) (\theta_{k}-\widehat{\theta }_{\mathcal
{A}%
,k})
\Biggr\rrvert,
\]
$\bolds{\theta}^{\ast}=t\widehat{\bolds{\theta
}}_{\mathcal
{A}%
} +(1-t)\bolds{\theta}$ for some $t\in{}[0,1]$.

Following \cite{qu:lindsay:li:00} and Lemma~S.4, for any
$\bolds{\theta}\in\Theta^{*}(C)$, one has
\begin{eqnarray*}
&&(\bolds{\theta}-\widehat{\bolds{\theta}}_{\mathcal
{A}})^{\T}
\dot{Q}%
_{n}(\widehat{\bolds{\theta}}_{\mathcal{A}})\\
&&\qquad=n(
\bolds{\theta }-\widehat{%
\bolds{\theta}}_{\mathcal{A}})^{\T}
\dot{\mathbf {G}}_{n}^{\T
}(\widehat{%
\bolds{\theta}}_{\mathcal{A}})\mathbf{C}_{n}^{-1}(
\widehat {\bolds{\theta}%
}_{\mathcal{A}})\mathbf{G}_{n}(
\widehat{\bolds{\theta }}_{\mathcal{A}%
}) \bigl\{ 1+o_{P}(1)
\bigr\}
\\
&&\qquad \leq C{ \rho_{n}^{-4}}n^{1/2}d_n
\bigl(n^{-1/2}+a_n\bigr),
\\
&&(\bolds{\theta}-\widehat{\bolds{\theta}}_{\mathcal
{A}})^{\T}
\ddot{Q}%
_{n}(\widehat{\bolds{\theta}}_{\mathcal{A}})
(\bolds{\theta }-\widehat{%
\bolds{\theta}}_{\mathcal{A}})
\\
&&\qquad =n(\bolds{\theta}-\widehat{\bolds{\theta}}_{\mathcal{A}%
})^{\T}
\dot{\mathbf{G}}_{n}^{\T}(\widehat{\bolds{\theta
}}_{\mathcal{A}})%
\mathbf{C}_{n}^{-1}(
\widehat{\bolds{\theta}}_{\mathcal
{A}})\dot{\mathbf{G}}%
_{n}(\widehat{\bolds{\theta}}_{\mathcal{A}}) (\bolds{
\theta }-\widehat{%
\bolds{\theta}}_{\mathcal{A}}) \bigl\{
1+o_{P}(1) \bigr\}
\\
&&\qquad  \geq C^2 { \rho_{n}^{-4}}
nd_n\bigl(n^{-1/2}+a_n\bigr)^2.
\end{eqnarray*}
By the Cauchy--Schwarz inequality,
$|R_{n}^{\ast}|\leq\frac{C^3}{6}  \{\sum_{i,j,k=1}^{d_n}
(\frac
{\partial^{3} Q_{n}(\bolds{\theta}^{*})}{\partial\theta_{i}\,
\partial\theta_{j} \,\partial\theta_{k}} )^2 \}^{1/2}
\times\|\bolds{\theta}-\widehat{\bolds{\theta
}}_{\mathcal{A}}\|^3$. According to assumption $n^{-1}d_{n}^{4}=o(1)$,
one has
%
\begin{eqnarray}\label{EQ:Qndiff}
\qquad &&Q_{n}(\bolds{\theta})-Q_{n}(\widehat{\bolds{
\theta }}_{\mathcal{A}})
\nonumber
\\[-8pt]
\\[-8pt]
\nonumber
&&\qquad= (\bolds{\theta}-\widehat{\bolds{\theta
}}_{\mathcal{A}})^{\T}\dot{Q}_{n}(\widehat{\bolds{
\theta }}_{\mathcal{A}})+\tfrac{1}{2}(\bolds{\theta}-\widehat {
\bolds{\theta}}_{\mathcal{A}})^{\T}\ddot{Q}_{n}(
\widehat{\bolds{\theta }}_{\mathcal{A}}) (\bolds{\theta}-\widehat{
\bolds{\theta }}_{\mathcal{A}}) \bigl\{ 1 + o_{P}(1) \bigr\}.
\end{eqnarray}
Thus, for sufficiently large $C$, the first term $(\bolds{\theta
}-\widehat{\bolds{\theta}}_{\mathcal{A}})^{\T}\dot
{Q}_{n}(\widehat
{\bolds{\theta}}_{\mathcal{A}})$ is dominated by the second term
$\frac{1}{2}(\bolds{\theta}-\widehat{\bolds{\theta
}}_{\mathcal{A}})^{\T}\ddot{Q}_{n}(\widehat{\bolds{\theta
}}_{\mathcal{A}})(\bolds{\theta}-\widehat{\bolds{\theta
}}_{\mathcal{A}})$.

Following the proof of Theorem~2 in \cite{xue:09}, if $\lambda
_{1n}\rightarrow0$, then for any $\bolds{\gamma}$ with $\Vert
\mathbf{B}^{%
\T} ( \bolds{\gamma}-\widehat{{\bolds{\gamma
}}}_{\mathcal{A}} )
\Vert_{n}=C_{2}{ \rho_{n}^{-3}}n^{-1/2}d_n^{1/2}$, one has $%
\Vert\bolds{\gamma}_{l}\Vert_{\mathbf{K}_{l}}\geq
a\lambda
_{1,n}$, and $%
\Vert\widehat{\bolds{\gamma}}_{\mathcal{A},l}\Vert
_{\mathbf
{K}_{l}}\geq
a\lambda_{1,n}$ for each $l=1,\ldots,s_{x}$, when $n$ is large enough.
By the definition of the SCAD penalty,
$\sum_{l=1}^{s_{x}} \{ p_{\lambda_{1,n}} ( \Vert
\bolds{\gamma}%
_{l}\Vert_{\mathbf{K}_{l}} ) -p_{\lambda_{1,n}} (
\Vert
\widehat{%
\bolds{\gamma}}_{\mathcal{A},l}\Vert_{\mathbf
{K}_{l}} )
\} =0$
for large $n$.\break  Furthermore, for any $\bolds{\beta}$ with $\Vert
\bolds{\beta}-\widehat{\bolds{\beta}}_{\mathcal{%
A}}\Vert\leq C{ \rho_{n}^{-3}}d_n^{1/2}(n^{-1/2}+a_{n})$,
\begin{eqnarray*}
&& \sum_{l=1}^{s_{z}}n\bigl\{
p_{\lambda_{2,n}} \bigl( |\beta _{l}| \bigr) -p_{\lambda_{2,n}}\bigl(|\widehat{
\beta}_{\mathcal{A},l}|\bigr)\bigr\}\\
&&\qquad=%
\sum_{l=1}^{s_{z}}np_{\lambda_{2,n}}^{\prime}\bigl(|
\widehat{\beta }_{\mathcal{A%
},l}|\bigr) (\beta_{l}-\widehat{
\beta}_{\mathcal{A},l}) \sgn(\widehat {\beta }_{%
\mathcal{A},l})
\\
&&\qquad\quad{} +\sum_{l=1}^{s_{z}}np_{\lambda_{2,n}}^{\prime\prime
}\bigl(|%
\widehat{\beta}_{\mathcal{A},l}|\bigr) ( \beta_{l}-\widehat{\beta
}_{\mathcal{A}%
,l}) ^{2}\sgn(\widehat{\beta}_{\mathcal{A},l})\bigl
\{1+o(1)\bigr\},
\end{eqnarray*}
$ 
p_{\lambda_{2,n}}^{\prime}(|\widehat{\beta}_{\mathcal
{A},l}|)-p_{\lambda
_{2n}}^{\prime}(|\beta_{0,l}|)=p_{\lambda_{2,n}}^{\prime\prime
}(|\beta
_{0,l}|)( \widehat{\beta}_{\mathcal{A},l}-\beta_{0,l}) \sgn(\beta
_{0,l})\{1+o(1)\}$. 
Thus,
\begin{eqnarray*}
&&\sum_{l=1}^{s_{z}}p_{\lambda_{2,n}}^{\prime
}\bigl(|
\widehat{\beta}_{0,l}|\bigr) ( \beta_{l}-\widehat{
\beta}_{\mathcal{A},l}) \sgn(\widehat{\beta }_{\mathcal{A
},l})\\
&&\qquad=\sum
_{l=1}^{s_{z}} \bigl\{ p_{\lambda_{2n}}^{\prime
}\bigl(|
\beta_{0,l}|\bigr) \bigr\} ( \beta_{l}-\widehat{
\beta}_{\mathcal{A},l}) \sgn( \widehat{\beta}_{\mathcal{A},l})
\\
&&\qquad \leq C{ \rho_{n}^{-3}} s_{z}^{1/2}a_nd_n^{1/2}
\bigl(n^{-1/2}+a_{n}\bigr) \leq C { \rho_{n}^{-3}}d_n
\bigl(n^{-1/2}+a_n\bigr)^2.
\end{eqnarray*}
Meanwhile,
\[
\sum_{l=1}^{s_{z}}np_{\lambda_{2,n}}^{\prime\prime}\bigl(|
\widehat {\beta} 
_{0,l}|\bigr) ( \beta_{l}-\widehat{
\beta}_{\mathcal{A},l}) ^{2}\sgn (\widehat {%
\beta}_{\mathcal{A},l}) \leq C^2 { \rho_{n}^{-6}}nb_n
d_n\bigl(n^{-1/2}+a_n\bigr)^2.
\]
Hence, $\sum_{l=1}^{s_{z}}n\{ p_{\lambda_{2,n}} ( |\beta
_{l}| )
-p_{\lambda_{2,n}}(|\widehat{\beta}_{\mathcal{A},l}|)\}$ is also dominated
by the second term of (\ref{EQ:Qndiff}). Hence, by choosing a
sufficiently large $C$, (\ref{equ_th2}) holds for large $n$. The proof
of Theorem~\ref{THM:SCAD-QIF} is completed.

\subsection{Proof of Theorem \texorpdfstring{\protect\ref{THM:sparsity}}{4}}
Let $\varrho_{n,d}=\rho_{n}^{-3}n^{-1/2}d_n^{1/2}$, and define
$\Theta_{1}=\{ \bolds{\theta}\dvtx \bolds{\theta}\in\Theta_{
\mathcal{A}},\Vert\bolds{\beta}-\bolds{\beta}_{0}\Vert=O_{P}(
\varrho_{n,d}),\Vert\mathbf{B}^{\T}( {\bolds{\gamma
}}-\widetilde{\bolds{\gamma}}) \Vert
_{n}=O_{P}(\varrho_{n,d})\}$, $\Theta_{l}=\{(\bolds{\beta
}_{1}^{\T
},\break \bolds{\bolds{\beta}}_{2}^{\T}, \bolds{\gamma
}^{\T
})^{\T}\dvtx  {\bolds{\beta}}_{1}={\bolds{\beta}}_{2}=%
\mathbf{0}, \bolds{\gamma}=(\mathbf{0},\ldots
,\mathbf
{0},{\bolds{\gamma}}%
_{l}^{\T},\mathbf{0},\ldots,\mathbf{0})^{\T}, \Vert
\mathbf
{B}^{\T}{\bolds{%
\gamma}}\Vert_{n}=O_{P}(\varrho_{n,d})\}$ for $l=s_{x}+1,\ldots,
d_{x}$ and
$\Theta_{l} = \{ (\bolds{\beta}_{1}^{\T},\bolds{\bolds{\beta}}_{2}^{\T%
},\bolds{\gamma}^{\T})^{\T}\dvtx  \bolds{\beta
}_{2}=\mathbf
{0}, {\bolds{\gamma
=0}}, \Vert\bolds{\beta}_{1}\Vert=O_{P}(\varrho
_{n,d})\}$
for $l=d_{x}+1$.
It suffices to show that uniformly for any $\bolds{\theta}\in
\Theta
_{1}$ and $\bolds{\theta}_{l}^{\ast}\in\Theta_{l}$,
$L_{n} (
\bolds{\theta} ) \leq L_{n} ( \bolds{\theta
}+\bolds{\theta}%
_{l}^{\ast} ) $, with probability $1$ as $n\rightarrow
\infty$ for any $s_{x}+1\leq l\leq d_{x}+1$. Observe that, for $%
l=s_{x}+1,\ldots,d_{x}$,
\begin{eqnarray*}
&&L_{n} \bigl( \bolds{\theta}+\bolds{\theta}_{l}^{\ast
}
\bigr) -L_{n}(%
\bolds{\theta}) 
\\
&&\qquad=Q_{n} \bigl( \bolds{\theta}+
\bolds{\theta}_{l}^{\ast
} \bigr) -Q_{n}(
\bolds{\theta})+np_{\lambda_{1,n}}^{\prime
} ( w_{l} ) \bigl(
\bigl\Vert\bolds{\gamma}^{\ast}_{l}\bigr\Vert _{\mathbf{K}%
_{l}}
\bigr)
\\
&&\qquad =\bolds{\gamma}_{l}^{\ast\T}\frac{\partial
}{\partial
\bolds{\gamma}_{l}}Q_{n}
(\bolds{\theta} ) +\frac
{1}{2}\bolds{%
\gamma}_{l}^{\ast\T}\frac{\partial}{\partial\bolds{\gamma
}_{l}\,\partial
\bolds{\gamma}_{l}^{\T}}Q_{n} (
\bolds{\theta} ) \bolds{\gamma}%
_{l}^{\ast}
\bigl\{ 1+o_{P}(1) \bigr\}
\\
&&\qquad\quad{} +np_{\lambda_{1,n}}^{\prime} ( w_{l} ) \bigl(\bigl \Vert
\bolds{\gamma}^{\ast}_{l}\bigr\Vert _{\mathbf{K}%
_{l}} \bigr)
\\
&&\qquad =n\lambda_{1,n}\bigl\Vert\mathbf{B}^{\T}{\bolds{
\gamma}} 
_{l}\bigr\Vert_{n} \biggl\{
\frac{R_{n}}{\lambda_{1,n}}+\frac{p_{\lambda
_{1n}}^{\prime} ( w_{l} ) }{\lambda_{1,n}} \biggr\} \bigl\{ 1+o_{P}(1)
\bigr\},
\end{eqnarray*}
where $w_{l}$ is a value between $0$ and $\Vert\bolds{\gamma}%
^{\ast}_{l}\Vert_{\mathbf{K}_{l}}$ and
\begin{eqnarray*}
R_{n}& = &n^{-1}\bigl\Vert\mathbf{B}^{\T}{
\bolds{%
\gamma}}_{l}\bigr\Vert_{n}^{-1}
\biggl\{\bolds{\gamma}_{l}^{\ast\T}\frac{\partial}{\partial
\bolds{%
\gamma}_{l}}Q_{n}
(\bolds{\theta} ) +\frac
{1}{2}\bolds{\gamma}%
_{l}^{\ast\T}\frac{\partial}{\partial\bolds{\gamma
}_{l}\,\partial\bolds{%
\gamma}_{l}^{\T}}Q_{n} (\bolds{
\theta} ) \bolds{\gamma}%
_{l}^{\ast} \bigl\{
1+o_{P}(1) \bigr\} \biggr\}
\\
& = &O_{P}\bigl({ \rho_n^{-4}}%
n^{-1/2}d_n^{1/2}\bigr).
\end{eqnarray*}
Noting that $R_{n}/\lambda_{1,n}=o_P(1)$, and $\liminf_{n\rightarrow
\infty}\lim\inf_{w\rightarrow0^{+}}p_{\lambda_{1n}}^{\prime} (
w ) /\lambda_{1,n}=1$, thus, uniformly for any $\bolds{\theta}\in
\Theta_{1}$ and $\bolds{\theta}_{l}^{\ast}\in\Theta_{l}$,
$L_{n} (
\bolds{\theta} ) \leq L_{n} ( \bolds{\theta
}+\bolds{\theta}%
_{l}^{\ast} ) $, with probability tending to $1$ as
$n\rightarrow
\infty$ for any $l=s_{x}+1,\ldots,d_{x}$. On the other hand, for $%
l=d_{x}+1$,
\begin{eqnarray*}
&&L_{n} \bigl( \bolds{\theta}+\bolds{\theta
}_{l}^{\ast} \bigr) -L_{n}(%
\bolds{
\theta}) 
\\
&&\qquad=\bolds{\beta}_{1}^{\ast\T}
\frac{\partial}{\partial
\bolds{\beta}_{1}}Q_{n} (\bolds{\theta} ) +\frac
{1}{2}
\bolds{%
\beta}_{1}^{\ast\T}\frac{\partial^{2}}{\partial\bolds{\beta
}%
_{1}\,\partial\bolds{\beta}_{1}^{\T}}Q (
\bolds{\theta } ) \bolds{%
\beta}_{1}^{\ast}
\bigl\{ 1+o_{P}(1) \bigr\}
\\
&&\qquad\quad{} +n\sum_{q=1}^{s_{z}}p_{\lambda_{2,n}}^{\prime}
\bigl( |w_{l,q}| \bigr) \beta_{q}^{\ast}\sgn\bigl(
\beta_{q}^{\ast}\bigr).
\end{eqnarray*}
Similar arguments show that uniformly for any $\bolds{\theta}\in
\Theta
_{1}$ and $\bolds{\theta}_{d_x+1}^{\ast}\in\Theta_{d_x+1}$, $L_{n}(
\bolds{\theta}) \leq L_{n}(\bolds{\theta}+\bolds{\theta}_{d_x+1}^{\ast
}) $, with probability tending to $1$ as $n\rightarrow\infty$. This
establishes the desired result. 

\subsection{Proof of Theorem \texorpdfstring{\protect\ref{THM:oracle}}{5}}

Let $\bolds{\beta}_{\mathcal{S}} = (\beta_{1},\ldots,\beta
_{s_{z}})^{\T}$,
$\bolds{\gamma}_{\mathcal{S}} = (\bolds{\gamma}_{1}^{\T
},\ldots,\bolds{%
\gamma}_{s_{x}}^{\T})^{\T}$. Denote $\bolds{\theta}_{\mathcal
{S}}=(%
\bolds{\beta}_{\mathcal{S}}^{\T},\bolds{\gamma
}_{\mathcal
{S}}^{\T})^{\T}$%
. In a similar way,\vspace*{1pt} define $\widetilde{\bolds{\theta
}}_{\mathcal
{S}0}=(%
\bolds{\beta}_{\mathcal{S}0}^{\T},\widetilde{\bolds{\gamma
}}_{\mathcal{S}%
}^{\T})^{\T}$, in which $\bolds{\beta}_{\mathcal{S}0}=(\beta
_{10},\ldots,\beta_{s_{z}0})^{\T}$ and $\widetilde{\bolds{\gamma
}}_{\mathcal{%
S}}=(\widetilde{\bolds{\gamma}}_{1}^{\T
},\ldots,\widetilde{\bolds{%
\gamma}}_{s_{x}}^{\T})^{\T}$. It can be shown easily that there
exist $
\widehat{\bolds{\beta}}_{\mathcal{S}}$ and $\widehat
{\bolds{\gamma}}_{%
\mathcal{S}}$ minimizing $L_{n} ( (\bolds{\beta}_{\mathcal
{S}}^{\T},%
\mathbf{0}^{\T})^{\T},(\bolds{\gamma}_{\mathcal{S}}^{\T
},\mathbf{0}^{\T})^{%
\T} ) $, that is,
\[
\frac{\partial}{\partial\bolds{\theta}_{\mathcal
{S}}}L_{n} ( \bolds{\theta} ) \bigg\rrvert
_{\bolds{\beta}=(
\widehat
{%
\bolds{\beta}}_{\mathcal{S}}^{\T},\mathbf{0}^{\T}) ^{\T
},
\bolds{%
\gamma}=( \widehat{\bolds{\gamma}}_{\mathcal{S}}^{\T
},\mathbf
{0}^{\T%
})^{\T}}=\mathbf{0.}
\]
In the following, we consider $L_{n} ( \cdot ) $ as a
function of $%
\bolds{\theta}_{\mathcal{S}}$, and denote $\dot{L}_{n}$ and
$\ddot{L}_{n}$
the gradient vector and Hessian matrix of $L_{n} ( \cdot )
$ with
respect to $\bolds{\theta}_{\mathcal{S}}$. The rest of the proof
follows similarly
as that of Theorem~\ref{THM:parametric}. Using Taylor expansion, one
has
\begin{eqnarray*}
\dot{L}_{n}\bigl(\widehat{\bolds{\theta}}_{\mathcal{S}}^{\pqif
}
\bigr)-\dot {L}_{n}(%
\widetilde{\bolds{
\theta}}_{\mathcal{S}0})& = &\ddot {L}_{n} \bigl(
\bolds{%
\theta}_{\mathcal{S}}^{\ast} \bigr) \bigl(
\widehat{\bolds{\theta }}_{\mathcal{S%
}}^{\pqif}-\widetilde{
\bolds{\theta}}_{\mathcal{S}0}\bigr)
\\
&&{} +\frac{1}{2}\bigl(\widehat{%
\bolds{
\theta}}_{\mathcal{S}}^{\pqif}-\widetilde{\bolds{
\theta}}_{%
\mathcal{S}0}\bigr)^{\T}\frac{\partial\dot{L}_{n} (
\bolds{\theta}_{%
\mathcal{S}} )}{\partial\bolds{\theta}_{\mathcal{S}}\,
\partial\bolds{%
\theta}_{\mathcal{S}}^{\T}}\bigg\rrvert
_{\bolds{\theta}_{\mathcal
{S}}=\bolds{%
\theta}_{\mathcal{S}}^{\ast}}\bigl(\widehat{\bolds{\theta }}_{\mathcal
{S}}^{%
\pqif}-
\widetilde{\bolds{\theta}}_{\mathcal{S}0}\bigr),
\end{eqnarray*}
where $\bolds{\theta}_{\mathcal{S}}^{\ast}=t\widehat
{\bolds{\theta}}_{%
\mathcal{S}}^{\pqif}+(1-t)\widetilde{\bolds{\theta
}}_{\mathcal
{S}0}$, for
some $t\in{}[0,1]$. Thus, we have
\begin{eqnarray*}
&&-n^{-1}\dot{Q}_{n}(\widetilde{\bolds{
\theta}}_{\mathcal{S}0})- \bolds{%
\kappa}_{n}(
\widetilde{\bolds{\theta}}_{\mathcal{S}0}) \\
&&\qquad=n^{-1}\bigl\{
\ddot{Q}%
_{n}(\bolds{\theta}_{\mathcal{S}0}) +
\bolds{\Lambda}( \bolds{\theta}_{%
\mathcal{S}0})\bigr\}\bigl(\widehat{
\bolds{\theta}}_{\mathcal
{S}}^{\pqif
}-\widetilde{%
\bolds{\theta}}_{\mathcal{S}0}\bigr)
\\
&&\qquad\quad{}+\frac{1}{2}n^{-1}\bigl(\widehat{\bolds{
\theta }}_{\mathcal{S}}^{%
\pqif}-\widetilde{\bolds{
\theta}}_{\mathcal{S}0}\bigr)^{\T}\frac
{\partial
\dot{Q}_{n} ( \bolds{\theta}_{\mathcal{S}}
)}{\partial
\bolds{%
\theta}_{\mathcal{S}}\, \partial\bolds{\theta}_{\mathcal
{S}}^{\T
}}\bigg\rrvert
_{%
\bolds{\theta}_{\mathcal{S}}=\bolds{\theta}_{\mathcal
{S}}^{\ast}}\bigl(%
\widehat{\bolds{\theta}}_{\mathcal{S}}^{\pqif}-
\widetilde {\bolds{\theta}}%
_{\mathcal{S}0}\bigr),
\end{eqnarray*}
where $\bolds{\kappa}_{n}^{\T}(\widetilde{\bolds{\theta
}}_{\mathcal{S}0})  =   ( \{p_{\lambda_{2,n}}^{\prime
}(|\beta_{0,k}|)\sgn(\beta
_{0,k}) \}_{k=1}^{s_z} ,  \{\frac{\partial}{\partial
\bolds{\gamma}_{l}}p_{\lambda_{1,n}} (
\Vert\bolds{\gamma}_{l}\Vert_{\mathbf{K}_{l}}
)|_{\bolds{\gamma}_{l}=\widetilde{\bolds{\gamma
}}_{l}}
\}_{l=1}^{s_{x}} )$\break
and $\bolds{\Lambda}(\bolds{\theta}_{\mathcal
{S}0})=\diag
( \{ p_{\lambda_{2,n}}^{\prime\prime}(|\beta_{0,k}|) \}
_{k=1}^{s_z},
\{\frac{\partial^2 }{\partial\bolds{\gamma}_{l}\,\partial
\bolds{\gamma}_{l}^{\T}}p_{\lambda_{1,n}} ( \Vert
\bolds{\gamma}_{l}\Vert
_{\mathbf{K}_{l}} )|_{\bolds{\gamma}_{l}=\widetilde
{\bolds{\gamma}}_{l}} \}_{l=1}^{s_{x}} )$.
Note that
\begin{eqnarray*}
&&\frac{\partial}{\partial\bolds{\gamma}_{l}}p_{\lambda
_{1,n}} \bigl( \Vert\bolds{\gamma}_{l}
\Vert_{\mathbf{K}_{l}} \bigr)  = p_{\lambda_{1,n}}^{\prime}\bigl
( \Vert\bolds{\gamma}_{l}\Vert_{\mathbf{K}_{l}} \bigr) \Vert\bolds{
\gamma}_{l}\Vert_{\mathbf{K}_{l}}^{-1}\mathbf
{K}_{l}\bolds{\gamma}_{l},
\\
&&\frac{\partial^2 }{\partial\bolds{\gamma}_{l}\,\partial
\bolds{\gamma}_{l}^{\T}}p_{\lambda_{1,n}}\bigl ( \Vert\bolds{\gamma }_{l}
\Vert _{\mathbf{K}_{l}} \bigr)\\
&&\qquad = p_{\lambda_{1,n}}^{\prime
} \bigl( \Vert \bolds{
\gamma}_{l}\Vert_{\mathbf{K}_{l}} \bigr) \Vert \bolds{
\gamma}_{l}\Vert_{\mathbf{K}_{l}}^{-1}\mathbf
{K}_{l}
\\
&&\qquad\quad{}+ \bigl\{ p_{\lambda_{1,n}}^{\prime\prime}\bigl ( \Vert \bolds{
\gamma}_{l}\Vert_{\mathbf{K}_{l}} \bigr) \Vert \bolds{
\gamma}_{l}\Vert_{\mathbf{K}_{l}}^{-2}-p_{\lambda
_{1,n}}^{\prime
}
\bigl( \Vert\bolds{\gamma}_{l}\Vert_{\mathbf
{K}_{l}} \bigr) \Vert
\bolds{\gamma}_{l}\Vert_{\mathbf{K}_{l}}^{-3} \bigr\}
\mathbf{K}_{l}\bolds{\gamma}_{l}\bolds{\gamma
}_{l}^{\T
}\mathbf{K}_{l},
\end{eqnarray*}
$J_{n}\rightarrow\infty$ and $\lambda_{1,n}\rightarrow0$, as
$n\rightarrow\infty$, so $\Vert\widetilde{\bolds{\gamma
}}_{l}\Vert_{\mathbf{K}_{l}}\geq a\lambda_{1,n}$ for $n$ large enough
and for each $l=1,\ldots,s_{x}$. Thus, $p_{\lambda_{1,n}}^{\prime} (
\Vert\widetilde{\bolds{\gamma}}_{l}\Vert_{\mathbf
{K}_{l}} ) =0$ and $p_{\lambda_{1,n}}^{\prime\prime} (
\Vert\widetilde{\bolds{\gamma}}_{l}\Vert_{\mathbf
{K}_{l}} ) =0$.

Similar to the proof of Theorem~\ref{THM:parametric}, one has
\begin{eqnarray*}
 \widehat{\bolds{\beta}}_{\mathcal
{S}}^{\pqif
} - \bolds{
\beta}_{\mathcal{S}%
0}&=& \bigl\{ \bigl(\mathbf{H}_{\bolds{\beta\beta}}(\bolds{
\theta}_{\mathcal{S}%
0}) + \bolds{\Lambda}_{\mathcal{S}0}\bigr) -
\mathbf {H}_{\bolds{\beta\gamma}}(%
\bolds{\theta}_{\mathcal{S}0})
\bigl(\mathbf {H}_{\bolds{\gamma\gamma}}(%
\bolds{
\theta}_{\mathcal{S}0}) \bigr)^{-1}\mathbf {H}_{\bolds{\gamma\beta}%
}(
\bolds{\theta}_{\mathcal{S}0}) \bigr\}^{-1}
\\
&&{} \times \bigl( \mathbf{I,H}_{\bolds{\gamma\beta}}\bigl(\bolds{%
\theta}_{\mathcal{S}}^{\ast}\bigr) \bigl( \mathbf{H}_{\bolds{\gamma\gamma}}
\bigl(%
\bolds{\theta}_{\mathcal{S}}^{\ast}\bigr) \bigr)
^{-1} \bigr) \bigl\{ \mathbf{S%
}_{n}(
\widetilde{\bolds{\theta}}_{\mathcal{S}0})+\bolds{
\kappa}_{n}(%
\widetilde{\bolds{\theta}}_{\mathcal{S}0})
\bigr\}
\\
&&{} + O_{P} \bigl(\rho_{n}^{-1}n^{-1}d_n
\bigr)+o_{P}\bigl(n^{-1/2}\bigr).
\end{eqnarray*}
Note that
\begin{eqnarray*}
&&\bigl(\mathbf{I,H}_{\bolds{\beta\gamma}}(\widetilde {\bolds{\theta
}}_{\mathcal{S}0})\mathbf{H}_{\bolds{\gamma\gamma
}}^{-1}(
\widetilde{%
\bolds{\theta}}_{\mathcal{S}0})\bigr) \bigl\{
\mathbf {S}_{n}(\widetilde {\bolds{%
\theta}}_{\mathcal{S}0})+\bolds{\kappa}_{n}(\widetilde {
\bolds{\theta}}%
_{0})\bigr\}
\\
&&\qquad=\bigl\{ \dot{\mathbf{G}}_{\bolds{\beta}}^{\T}(\widetilde {
\bolds{%
\theta}}_{\mathcal{S}0})-\mathbf{H}_{\bolds{\beta\gamma
}}(
\widetilde{%
\bolds{\theta}}_{\mathcal{S}0})\mathbf{H}_{\bolds{\gamma
\gamma}}^{-1}(%
\widetilde{\bolds{\theta}}_{\mathcal{S}0})\dot{\mathbf
{G}}_{\bolds{\gamma
}}^{\T}(\widetilde{\bolds{\theta}}_{\mathcal{S}0})
\bigr\} \mathbf{C}_{n}^{-1}(\widetilde{\bolds{
\theta}}_{\mathcal
{S}0})\mathbf{G}%
_{n}(\widetilde{
\bolds{\theta}}_{\mathcal{S}0})+\bolds{\kappa}_{\mathcal{%
S}}.
\end{eqnarray*}
The asymptotic distribution of $\sqrt{n}\mathbf{A}_{n}\bolds{\Sigma}_{%
\mathcal{S},n}^{-1/2}(\widehat{\bolds{\beta}}_{\mathcal
{S}}^{\pqif
}-%
\bolds{\beta}_{\mathcal{S}0})$ is the same as that of
\begin{eqnarray*}
&&\sqrt{n}\mathbf{A}_{n}\bolds{\Sigma}_{\mathcal
{S},n}^{-1/2}
\bigl\{\mathbf{H}_{%
\bolds{\beta\beta}}(\widetilde{\bolds{\theta
}}_{\mathcal
{S}0})-\mathbf{H%
}_{\bolds{\beta\gamma}}(\widetilde{
\bolds{\theta }}_{\mathcal{S}0})%
\mathbf{H}_{\bolds{\gamma\gamma}}^{-1}(
\widetilde {\bolds{\theta}}_{%
\mathcal{S}0})\mathbf{H}_{\bolds{\gamma\beta}}(
\widetilde {\bolds{\theta}%
}_{\mathcal{S}0})+\bolds{
\Lambda}_{\mathcal{S}0}\bigr\}^{-1}
\\
&&\qquad{} \times \bigl\{ \bigl( \dot{\mathbf {G}}_{\bolds{\beta}}^{\T}(%
\widetilde{\bolds{\theta}}_{\mathcal{S}0}) - \mathbf
{H}_{\bolds{\beta
\gamma}}(\widetilde{\bolds{\theta}}_{\mathcal
{S}0})\mathbf
{H}_{\bolds{%
\gamma\gamma}}^{-1}(\widetilde{\bolds{\theta}}_{\mathcal
{S}0})
(%
\widetilde{\bolds{\theta}}_{\mathcal{S}0})\dot{\mathbf
{G}}_{\bolds{\gamma
}}^{\T}(\widetilde{\bolds{\theta}}_{\mathcal{S}0})
\bigr) \mathbf{C}%
_{n}^{-1}(\widetilde{
\bolds{\theta}}_{\mathcal
{S}0})\mathbf {G}_{n}(%
\widetilde{\bolds{\theta}}_{\mathcal{S}0}) \\
&&\hspace*{310pt}{}+ \bolds{\kappa
}_{\mathcal{S}%
} \bigr\}.
\end{eqnarray*}
Next, write $\widehat{\mathbf{J}}_\DZcalS^{\T}=\{(\widehat
{\mathbf{J}}_\DZcalS%
^{(1)})^{\T},\ldots,(\widehat{\mathbf{J}}_\DZcalS^{(K)})^{\T}\}
^{\T
}$, where $%
\widehat{\mathbf{J}}_\DZcalS^{(k)}=\frac{1}{n}\sum_{i=1}^{n}\mathbf{D}_{i}^{%
\T}\bolds{\Gamma}_{0,i}^{(k)} \widehat{\mathbf
{Z}}_{\mathcal
{S}i}$ and $%
\widehat{\mathbf{Z}}_{\mathcal{S}i}=\mathbf{Z}_{\mathcal
{S}i}-\mathbf{B}_{%
\mathcal{S}i}\{ \mathbf{J}_{\DBcalS}^{\T}( \mathbf{C}_{n}^{0})
^{-1}\mathbf{%
J}_{\DBcalS}\} ^{-1}\mathbf{J}_{\DBcalS}^{\T}(\mathbf
{C}_{n}^{0}) ^{-1}%
\mathbf{J}_{\DZcalS}$. Then we can express
\[
\mathbf{H}_{\bolds{\beta\beta}}(\widetilde{\bolds{\theta}}_{\mathcal{S}%
0})-
\mathbf{H}_{\bolds{\beta\gamma}}(\widetilde {\bolds{\theta}}_{%
\mathcal{S}0})
\mathbf{H}_{\bolds{\gamma\gamma
}}^{-1}(\widetilde{\bolds{%
\theta}}_{\mathcal{S}0})\mathbf{H}_{\bolds{\gamma\beta
}}(\widetilde{%
\bolds{\theta}}_{\mathcal{S}0})=\widehat{\mathbf {J}}_{\DZcalS
}^{\T}
\bigl( \mathbf{C}_{n}^{0} \bigr) ^{-1}
\widehat{\mathbf{J}}_\DZcalS ^{\T%
}\bigl
\{1+o_{P}(1)\bigr\}.
\]
Using similar arguments as given in Lemma~S.11, we know
\begin{eqnarray*}
&&\mathbf{A}_{n}\bolds{\Sigma}_{\mathcal{S},n}^{-1/2}
\bigl\{ \mathbf{H}_{\bolds{\beta\beta}}(\widetilde{\bolds{
\theta}}_{\mathcal{S}0})-\mathbf{H}_{\bolds{\beta\gamma
}}(\widetilde{\bolds{
\theta }}_{\mathcal{S}0})\mathbf{H}_{\bolds{\gamma\gamma
}}^{-1}(
\widetilde{\bolds{\theta}}_{\mathcal{S}0})\mathbf {H}_{\bolds{\gamma\beta}}(
\widetilde{\bolds{\theta }}_{\mathcal{S}0})+\bolds{
\Lambda}_{\mathcal{S}0}\bigr\}^{-1}
\\
&&\quad{} \times \bigl( \dot{\mathbf{G}}_{\bolds{\beta}}^{\T}(\widetilde {
\bolds{\theta}}_{\mathcal{S}0})-\mathbf{H}_{\bolds{\beta\gamma}}(\widetilde{
\bolds{\theta }}_{\mathcal{S}0})\mathbf{H}_{\bolds{\gamma\gamma
}}^{-1}(
\widetilde{\bolds{\theta}}_{\mathcal{S}0})\dot{\mathbf
{G}}_{\bolds{\gamma}}^{\T}(\widetilde{\bolds{\theta
}}_{\mathcal{S}0}) \bigr) \mathbf{C}_{n}^{-1}(
\widetilde {\bolds{\theta}}_{\mathcal{S}0})\mathbf{G}_{n}(
\widetilde{\bolds{\theta}}_{\mathcal{S}0})
\\
&& \qquad=\mathbf{A}_{n}\bolds{\Sigma}_{\mathcal
{S},n}^{-1/2}
\bigl\{ \mathbf{H}_{\bolds{\beta\beta}}(\widetilde{\bolds{
\theta}}_{\mathcal{S}0})-\mathbf{H}_{\bolds{\beta\gamma
}}(\widetilde{\bolds{
\theta }}_{\mathcal{S}0})\mathbf{H}_{\bolds{\gamma\gamma
}}^{-1}(
\widetilde{\bolds{\theta}}_{\mathcal{S}0})\mathbf {H}_{\bolds{\gamma\beta}}(
\widetilde{\bolds{\theta }}_{\mathcal{S}0})+\bolds{
\Lambda}_{\mathcal{S}0}\bigr\}^{-1}
\\
&& \qquad\quad{}\times\widehat{\mathbf{J}}_{\DZcalS}^{\T} \bigl(
\mathbf{C}_{n}^{0} \bigr) ^{-1}
\mathbf{G}_{n}^{0} +o_{P}
\bigl(n^{-1/2}\bigr).
\end{eqnarray*}
Thus, the desired result follows. 
\end{appendix}

\section*{Acknowledgements}
The authors thank the Co-Editor, an Associate Editor and three referees
for their constructive comments that have substantially improved an
earlier version of this paper.
L. Wang and L. Xue contribute equally for the paper.


\begin{supplement}[id=suppA]
\stitle{Supplement to ``Estimation and model selection in generalized additive partial
linear models for correlated data with diverging number of covariates''}
\slink[doi]{10.1214/13-AOS1194SUPP} 
\sdatatype{.pdf}
\sfilename{aos1194\_supp.pdf}
\sdescription{The supplementary material provides a number of technical
lemmas and their proofs. The technical lemmas are used in the proofs of
Theorems \ref{THM:nonparametric}--\ref{THM:oracle} in the paper.}
\end{supplement}



%

\printaddresses

\end{document}